\documentclass[12pt,centertags,oneside]{amsart}

\usepackage{amsmath,amstext,amsthm,amscd,typearea,hyperref,stmaryrd,mathabx}
\usepackage{amssymb}
\usepackage{a4wide}
\usepackage[mathscr]{eucal}
\usepackage{mathrsfs}
\usepackage{typearea}
\usepackage{charter}
\usepackage{pdfsync}
\usepackage[a4paper,width=16.2cm,top=3cm,bottom=3cm,lines=50]{geometry}
\usepackage[all]{xy}
\usepackage{tikz}

\numberwithin{equation}{section}

\title[]{Poincar\'e metric of holomorphic foliations with non-degenerate singularities}

\author{Fran\c cois Bacher}
\address{Universit\'e de Lille, 
Laboratoire de math\'ematiques Paul Painlev\'e, 
CNRS U.M.R. 8524,  
59655 Villeneuve d'Ascq Cedex, 
France.}
\email{francois.bacher@univ-lille.fr}

\date{\today}

\theoremstyle{plain}
\newtheorem{thm}{Theorem}[section]
\newtheorem{lem}[thm]{Lemma}

\newtheorem{prop}[thm]{Proposition}

\newtheorem*{thm*}{Theorem}
\newtheorem*{conj*}{Conjecture}

\theoremstyle{definition}
\newtheorem{defn}[thm]{Definition}
\newtheorem{exmp}[thm]{Example}
\newtheorem*{exmp*}{Example}

\theoremstyle{remark}
\newtheorem{rem}[thm]{Remark}

\newtheorem*{pr}{Proof}

\makeatletter

\@addtoreset{equation}{section}
\makeatother

\DeclareMathOperator{\id}{id}

\DeclareMathOperator{\Jac}{Jac}

\DeclareMathOperator{\Tr}{Tr}

\DeclareMathOperator{\com}{com}

\newcommand{\cjg}[1]{\overline{#1}}

\newcommand{\adh}[1]{\overline{#1}}

\newcommand{\eps}{\varepsilon}

\newcommand{\fol}{\mathscr{F}}

\newcommand{\set}[1]{\mathbb{#1}}

\newcommand{\germ}[2]{\left(#1,#2\right)}

\newcommand{\der}[2]{\frac{\partial#1}{\partial#2}}

\newcommand{\disder}[2]{\displaystyle\der{#1}{#2}}

\newcommand{\bder}[2]{\der{#1}{\cjg{#2}}}

\newcommand{\disbder}[2]{\displaystyle\bder{#1}{#2}}

\newcommand{\deder}[3]{\frac{\partial^2#1}{\partial#2\partial#3}}

\newcommand{\dbder}[3]{\deder{#1}{#2}{\cjg{#3}}}

\newcommand{\dder}[2]{\frac{\partial^2#1}{\partial#2^2}}

\newcommand{\bdder}[2]{\dder{#1}{\cjg{#2}}}

\newcommand{\leafatlas}{\mathscr{L}}

\newcommand{\leaf}{L}

\newcommand{\leafu}[1]{\leaf_{#1}}

\newcommand{\ballleafu}[2]{\leafu{#1}{\left[#2\right]}}

\newcommand{\norm}[1]{\left\Vert#1\right\Vert}

\newcommand{\inter}[2]{\left[#1,#2\right]}

\newcommand{\intero}[2]{\left(#1,#2\right)}

\newcommand{\intlo}[2]{\left(#1,#2\right]}

\newcommand{\Cmod}[1]{\left\vert#1\right\vert}

\newcommand{\hol}{\mathcal{H}}

\newcommand{\holo}[2]{\hol\left(#1,#2\right)}

\newcommand{\projective}[1]{\set{P}^{#1}}

\newcommand{\class}[1]{\mathscr{C}^{#1}}

\newcommand{\smooth}{\class{\infty}}

\newcommand{\flot}[1]{\varphi_{#1}}

\newcommand{\mani}[1]{#1}

\newcommand{\manis}[2]{\mani{#1}\backslash\mani{#2}}

\newcommand{\PC}{P}

\newcommand{\dsing}[2]{d\left(#2,\mani{#1}\right)}

\newcommand{\rD}[1]{#1\set{D}}

\newcommand{\adhrD}[1]{#1\adh{\set{D}}}

\newcommand{\DR}[1]{\set{D}_{#1}}

\newcommand{\adhDR}[1]{\adh{\set{D}}_{#1}}

\newcommand{\dherm}{d}

\newcommand{\dhermset}[1]{d_{#1}}

\newcommand{\dPC}{d_{\PC}}

\newcommand{\dPCset}[1]{d_{\PC,#1}}

\newcommand{\dhermf}[1]{d_{\leafu{#1}}}

\newcommand{\dhermfset}[2]{d_{\leafu{#1},#2}}

\newcommand{\setst}{~\vert~}

\newcommand*{\transp}[1]{\left(#1\right)^T}

\begin{document}

\theoremstyle{plain}

\begin{abstract} Consider a Brody hyperbolic foliation $\fol{}$ with non-degenerate singularities on a compact complex manifold. We show that its leafwise Poincar\'e metric is transversally H\"{o}lder continuous with a logarithmic slope towards the singular set of $\fol{}$.

  \end{abstract}

\maketitle

\section{Introduction}

The ergodic theory of laminations by Riemann surfaces has received a lot of attention during the last two decades. Much progress has been focused on the global dynamical properties of laminations by hyperbolic Riemann surfaces with singularities, and in particular singular holomorphic foliations. Indeed, most of interesting holomorphic foliations are singular and hyperbolic. In the case of $\projective{n}$, all foliations are singular. For $d\in\set{N}$, let $\fol_d\left(\projective{n}\right)$ be the space of foliations of degree $d$ on $\projective{n}$. By the results of Lins Neto and Soares~\cite{LNS} and Jouanolou~\cite{Jou}, the singularities of a generic foliation $\fol{}\in\fol_d\left(\projective{n}\right)$, for $d\geq2$, are all non-degenerate. Moreover, by the results of Glutsyuk~\cite{Glu} and Lins Neto~\cite{LN2}, if all the singularities of a foliation $\fol{}\in\fol_d\left(\projective{n}\right)$, for $d\geq2$, are non-degenerate, then the foliation $\fol{}$ is hyperbolic. It is even Brody hyperbolic in the sense of~\cite{DNSII}. Moreover, Loray and Rebelo build in~\cite{LorReb} a non-empty open set in $\fol_d{\left(\projective{n}\right)}$ in which every foliation has all its leaves that are dense. When $n=2$, Nguyên computes the Lyapunov exponent of a generic foliation $\fol\in\fol_d\left(\projective{2}\right)$ in \cite{NguLyap1,NguLyap2}. We recall briefly some recent developments and refer the reader to the survey articles \cite{surForSib,surDinhSib,surVANG18,surVANG21} for more detailed expositions.

By solving heat equations with respect to a positive harmonic current, Dinh, Nguyên and Sibony obtain in~\cite{DNS12} abstract geometric ergodic theorems for laminations and singular holomorphic foliations. As a consequence, they obtain a concrete and effective Birkhoff ergodic theorem in the case of holomorphic foliations with linearizable singularities. In a series of two articles~\cite{DNSI,DNSII}, these authors study the modulus of continuity of the leafwise Poincar\'e metric for a compact lamination and for a holomorphic foliation with linearizable singularities. Using this transversal regularity, they introduce the hyperbolic entropy and establish the finiteness of this entropy for a large class of foliations in dimension 2 with linearizable singularities.

Let us recall more precisely their results about the Poincar\'e metric. Let $M$ be a compact complex manifold and $\fol{}=\left(M,\leafatlas{},E\right)$ be a smooth foliation by hyperbolic Riemann surfaces. That is, $\leafatlas{}$ is an atlas of flow boxes of $M\backslash E$ and $E$ is the singular set of $\fol{}$. Let $g_{\PC}$ be the leafwise Poincar\'e metric of $\fol{}$ and fix an ambient smooth Hermitian metric $g_{\mani{M}}$ on $\mani{M}$. Consider then the function $\eta\colon M\backslash E\to\intero{0}{+\infty}$ given by $\eta^2g_{\PC}=4g_{\mani{M}}$, where on the right hand side $g_{\mani{M}}$ is restricted to each leaf of $\fol{}$. The normalization of this function will be explained in Subsection~\ref{secPCfactor}. The two main theorems of Dinh, Nguy\^{e}n and Sibony we are interested in are the following.

\begin{thm}[Dinh-Nguy\^{e}n-Sibony~\cite{DNSI}] \label{thmDNSI}Let $\fol{}=\left(M,\leafatlas{}\right)$ be a transversally smooth compact foliation by hyperbolic Riemann surfaces. Then, the function $\eta$ is transversally H\"{o}lder continuous. Moreover, the H\"{o}lder exponent can be expressed in geometric terms.
\end{thm}

\begin{thm}[Dinh-Nguy\^{e}n-Sibony~\cite{DNSII}] \label{thmDNSII} Let $\fol{}=\left(M,\leafatlas{},E\right)$ be a Brody hyperbolic singular holomorphic foliation on a compact complex manifold $\mani{M}$. Suppose that all singularities of $\fol{}$ are linearizable. Then, there exist constants $C>0$ and $\alpha\in\intero{0}{1}$ such that
  \[\Cmod{\eta(x)-\eta(y)}\leq C\left(\frac{\max\left(\log^{\star}\dsing{E}{x},\log^{\star}\dsing{E}{y}\right)}{\log^{\star}d(x,y)}\right)^{\alpha},\quad x,y\in\manis{M}{E},\]
  where $\log^{\star}=1+\Cmod{\log}$ is a log-type function, $\dsing{E}{x}$ denotes the distance from $x$ to the singular set $\mani{E}$ and $d(\cdot,\cdot)$ is the distance induced by the ambient Hermitian metric $g_{\mani{M}}$.
\end{thm}

So, it is natural to wonder whether these basic questions can be solved when the singularities of $\fol{}$ are more general. The purpose of this work is to provide a positive answer to the problem of transversal regularity of the Poincar\'e metric when the singularities are all non-degenerate. Our main result is the following.

\begin{thm} \label{mainthm} Let $\fol{}=\left(M,\leafatlas{},E\right)$ be a Brody hyperbolic singular holomorphic foliation on a compact complex manifold $\mani{M}$ with Hermitian metric $g_{\mani{M}}$. Suppose that all the singularities of $\fol{}$ are non-degenerate. Then, there exist constants $C>0$ and $\alpha\in\intero{0}{1}$ such that
  \begin{equation}\label{eqmainthm}\Cmod{\eta(x)-\eta(y)}\leq C\left(\frac{\max\left(\log^{\star}\dsing{E}{x},\log^{\star}\dsing{E}{y}\right)}{\log^{\star}d(x,y)}\right)^{\alpha},\quad x,y\in\manis{M}{E},\end{equation}
  where we use the same notations as in Theorem~\ref{thmDNSII}. 
\end{thm}
  
Let us explain the method of our proof. We basically follow the approach of Dinh, Nguy\^{e}n and Sibony in~\cite{DNSII} using the Beltrami equation. However, their proof uses the peculiar geometry of a linearizable singularity as a local model in order to study the regularity of the Beltrami coefficient near the singularities. Namely, they use a sort of local invariance by homothety of a linearizable singularity. In the context of a non-degenerate singularity, such a local model is not available anymore. The novelty of our work is to settle the same kind of estimates as theirs in this more general situation.

The article is organized as follows. In Section~\ref{secprelim}, we recall some definitions about hyperbolic singular holomorphic foliations and the Poincar\'{e} metric. We specify the type of singularities we will consider, to point out our improvements to Theorem~\ref{thmDNSII}. We also state some general results, which will be useful in our proof. Next, we build in Section~\ref{secorthproj} a local orthogonal projection from a leaf onto another in singular charts. We get precise estimates on its $\class{2}-$norm using that the singularities are non-degenerate. More precisely, we use only the speed of expansion of a flow in a leaf $\leafu{x}$ with respect to the Hermitian metric $g_{\mani{M}}$. We get the same kind of estimates as in the case of linearizable singularities. Finally, we prove Theorem~\ref{mainthm} in Section~\ref{secmainthm}. We follow closely the proof of Theorem~\ref{thmDNSII}, with a slightly different exposition.

\subsection*{Notations} Throughout this paper, we will denote by $\set{D}$ the unit disk of $\set{C}$, and $\rD{r}$ (respectively $\adhrD{r}$) the open (respectively closed) disk of radius $r\in\set{R}_+^*$ for the standard Euclidean metric of $\set{C}$. For $R\in\set{R}_+^*$, we will also denote by $\DR{R}$ (respectively $\adhDR{R}$) the open (respectively closed) disk of hyperbolic radius $R$ in $\set{D}$, so that $\DR{R}=\rD{r}$ with $r=\frac{e^R-1}{e^R+1}$, or if $r\in\intero{0}{1}$, with $R=\ln\frac{1+r}{1-r}$. More generally, for $\rho\in\set{R}_+^*$ and $U$ a subset of a vector space, $\rho U$ will denote the image of $U$ by the homothety $z\mapsto\rho z$.

Throughout this paper, we will denote by $C$, $C'$, $C''$, etc\dots{} positive constants which will not always be the same.

\subsection*{Acknowledgments} The  author is supported by the Labex CEMPI (ANR-11-LABX-0007-01) and by the project QuaSiDy (ANR-21-CE40-0016). We would like to thank Julio Rebelo for helping us to find Example~\ref{exfin}. We are grateful to the referee for his helpful remarks.

\section{Preliminaries \label{secprelim}}

\subsection{Poincar\'e metric of a holomorphic foliation\label{secPCfactor}}

Let $\fol{}=\left(M,\leafatlas{},E\right)$ be a singular holomorphic foliation by curves on a complex manifold $\mani{M}$. For $x\in\manis{M}{E}$, we denote by $\leafu{x}$ the leaf of $\fol{}$ passing through $x$. 

For two complex manifolds $\mani{N}$ and $\mani{N}'$, let $\holo{\mani{N}}{\mani{N}'}$ be the set of holomorphic functions from $\mani{N}$ to $\mani{N}'$. Define
\[\holo{\set{D}}{\fol{}}=\left\{f\in\holo{\set{D}}{\manis{M}{E}}\setst{}f(\set{D})\subset\leafu{f(0)}\right\}.\]

Fix a Hermitian metric $g_{\mani{M}}$ on $\mani{M}$ and define the function $\eta\colon\manis{M}{E}\to\intlo{0}{+\infty}$ as follows.
\[\eta(x)=\sup\left\{\norm{\alpha'(0)}_{g_{\mani{M}}}\setst{}\alpha\in\holo{\set{D}}{\fol{}},\,\alpha(0)=x\right\},\]
% \begin{equation}\label{deflambda}\eta(x)=\sup\left\{\norm{\alpha'(0)}_{g_{\mani{M}}}\setst{}\alpha\in\holo{\set{D}}{\fol{}},\,\alpha(0)=x\right\},\end{equation}
where $\norm{v}_{g_{\mani{M}}}$ is the norm of $v\in T_x\mani{M}$ with respect to the metric $g_{\mani{M}}$, that is $\norm{v}_{g_{\mani{M}}}=\sqrt{g_{\mani{M},x}(v,v)}$.

The function $\eta$ is defined to satisfy the following facts, proven by Verjovsky~\cite{Ver}.

\begin{prop}\label{impdeflambda}
  \begin{enumerate}
  \item\label{ifflfinite} For $x\in\manis{M}{E}$, $\eta(x)<+\infty$ if and only if the leaf $\leafu{x}$ is hyperbolic, that is, it is uniformized by the Poincar\'e disk $\set{D}$.

   \item\label{valuelambda} If $\leafu{x}$ is hyperbolic, we have $\eta(x)=\norm{u'(0)}_{g_{\mani{M}}}$, where $u\colon\set{D}\to\leafu{x}$ is any uniformization of $\leafu{x}$ such that $u(0)=x$.
  \item\label{genPoincaremetric} If $\leafu{x}$ is hyperbolic, then $\frac{4g_{\mani{M}}}{\eta^2}$ induces the Poincar\'e metric on $\leafu{x}$.
  \end{enumerate}
\end{prop}

In the general case, the regularity of the function $\eta$ is very weak. In fact, it is not always continuous. Forn\ae{}ss and Sibony have only shown in \cite[Theroem 20]{surForSib} that it is lower semicontinuous if $\mani{M}$ is compact. They also give sufficient conditions for $\eta$ to be continuous. In our context, they are not necessarilly satisfied.

We recall from~\cite{DNSII} the following.

\begin{defn} \label{defBrody}

We say that $\fol{}$ is \emph{Brody hyperbolic} if there exists a positive constant $A$ such that $\eta<A$ on $\manis{M}{E}$.

\end{defn}

\begin{defn}
  Near a singularity $p$, there exists an open chart $U\subset\mani{M}$ with coordinates $z\in\set{C}^n$, such that the leaves of $\fol{}$ are locally defined as the complex flow of a holomorphic vector field
  \[X=\sum_{j=1}^nF_j\der{}{z_j},\qquad F_j\in\holo{U}{\set{C}},~1\leq j\leq n.\]

  If in some chart centered at $p$, $X$ is of the form $\sum_{j=1}^n\lambda_jz_j\der{}{z_j}$ with all $\lambda_j\neq0$, we say that $p$ is a \emph{linearizable singularity}.

  The functions $F_j$ can be developed as a power series $F_j=\sum_{\alpha\in\set{N}^n}c_{\alpha,j}z^{\alpha}$. The \emph{1-jet} of $X$ at $p$ is defined in the chart $(U,z)$ as
  \[X_1=\sum\limits_{j=1}^n\sum\limits_{\Cmod{\alpha}\leq1}c_{\alpha,j}z^{\alpha}\der{}{z_j}.\]

  See for example~\cite[Chapter I]{Ilya} for more details. If the 1-jet of $X$ has an isolated singularity at $p$, we say that $p$ is a \emph{non-degenerate singularity}.
\end{defn}

\subsection{Some useful results}

Let us prove the following lemma, which gives a bound to the Poincar\'e distance on a cut-off disk if the Euclidean distance is very small. For $z,w\in\set{D}$, we denote by $\dPC{}(z,w)$ the Poincar\'e distance between $z$ and $w$.

\begin{lem}\label{lemPCeucl} Let $z,w\in\DR{R}$ with $R\in\set{R}_+^*$. If $\Cmod{z-w}\leq e^{-R}$ and if $R$ is sufficiently large, then $\dPC{}(z,w)\leq e^R\Cmod{z-w}$.

\end{lem}

\begin{pr} Take $r=\frac{e^R-1}{e^R+1}$ the corresponding Euclidean radius such that $\Cmod{z}\leq r$ and $\Cmod{w}\leq r$. Using the straight line from $z$ to $w$, we get the following bound for the Poincar\'e distance between $z$ and $w$.
  \[\dPC{}(z,w)\leq\frac{2\Cmod{z-w}}{1-r^2}\leq e^R\Cmod{z-w},\]
if $R$ is sufficiently large, since $1-r^2=\frac{4e^R}{\left(e^R+1\right)^2}$.
  \qed
\end{pr}

We also recall the following version of the Gr\"{o}nwall Lemma, which will be useful for our further work. See for example~\cite{MP} for a proof. For $f\colon U\to V$ a $\class{1}$ function between two vector spaces or manifolds, and for $u\in U$, we denote by $d_uf$ the differential of $f$ with respect to $u$.

  \begin{lem}[Gr\"{o}nwall] \label{CGron}
    Let $\varphi$ be a $\class{1}$ function from an open ball $B_r$ centered at the origin of radius $r$ in $\set{R}^m$ to $\set{C}^n$. Let also $C$ be a non negative constant. If we suppose that for all $x\in B_r$,
    \[\norm{d_x\varphi}\leq C\norm{\varphi(x)},\]
    then for $x\in B_r$, we have
    \[\norm{\varphi(x)}\leq\norm{\varphi(0)}e^{C\norm{x}}.\]
  \end{lem}

  \section{Local study of the orthogonal projection onto a leaf \label{secorthproj}}

  We want to study the orthogonal projection from a leaf $\leafu{x}$ near $x$ onto another leaf $\leafu{y}$ near $y$ inside a singular chart of $\fol{}$. Consider then a foliation $\fol{}$ on $3\rho\set{D}^n$ for some $\rho>0$ with a non-degenerate singularity at the origin. We suppose that $\fol{}$ is generated by a holomorphic vector field $X$ on $3\rho\set{D}^n$. We endow $3\rho\set{D}^n$ with the usual Hermitian product $\left\langle\cdot,\cdot\right\rangle$ of $\set{C}^n$ and the associated norm $\norm{\cdot}$. Since $0$ is a non-degenerate singularity, there exists a constant $C_0>1$ such that
  \begin{equation}\label{borneX}C_0^{-1}\norm{z}\leq\norm{X(z)}\leq C_0\norm{z},\quad z\in3\rho\set{D}^n.\end{equation}

  For $z\in2\rho\set{D}^n$, denote by $\flot{z}$ the flow of $X$ starting at $z$. Namely, $\flot{z}\colon\germ{\set{C}}{0}\to3\rho\set{D}^n$ is a maximal solution of the Cauchy problem
  %\begin{equation}\label{defCauchypb}\left\{\begin{aligned}\der{\flot{z}}{t}(t)&=X(\flot{z}(t))\\\flot{z}(0)&=z.\end{aligned}\right.\end{equation}
  \[\left\{\begin{aligned}\der{\flot{z}}{t}(t)&=X(\flot{z}(t))\\\flot{z}(0)&=z.\end{aligned}\right.\]
  
  Taking the derivative of the first line with respect to $t$, we find that there exists a constant $C_1>0$ such that
  \begin{equation}\label{bornederflot} \norm{\dder{\flot{z}}{t}(t)}\leq C_1\norm{X(\flot{z}(t))},\quad\norm{\frac{\partial^3\flot{z}}{\partial^3t}(t)}\leq C_1\norm{X(\flot{z}(t))}.\end{equation}

  Since we have $\norm{\der{\flot{z}}{t}(t)}\leq C_0\norm{\flot{z}(t)}$, it follows from Gr\"{o}nwall Lemma that $\flot{z}$ is at least defined on a disk of radius $r_0=\frac{1}{C_0}\ln\frac{3}{2}$. Again by Gr\"{o}nwall Lemma, it is also clear that for $t\in\rD{r_0}$, $\norm{\flot{z}(t)}\geq\frac{2}{3}\norm{z}$. In particular, there exists a constant $C>1$ such that $C^{-1}\norm{X(z)}\leq\norm{X(\flot{z}(t))}\leq C\norm{X(z)}$ for $t\in\rD{r_0}$.

  We first show that the local orthogonal projection exists, is unique, and that it is $\class{1}$. For $x,y\in 2\rho\set{D}^n\backslash\{0\}$, define
  \begin{equation}\label{defgxy}g(t,u)=\bder{\norm{\flot{x}(t)-\flot{y}(u)}^2}{u}(t,u)=\left\langle \flot{x}(t)-\flot{y}(u),X(\flot{y}(u))\right\rangle.\end{equation}
  For $t,u\in r_0\set{D}$, if $g(t,u)=0$, then the vector $\flot{x}(t)-\flot{y}(u)$ is orthogonal to the tangent line to $\leafu{y}$ in $\flot{y}(u)$. That is how we will characterize the local orthogonal projection.

  \begin{lem} \label{existf} There exist a radius $r_1<r_0$ and a positive number $\eps_0$ satisfying the following conditions. Let $x,y$ be two points in $2\rho\set{D}^n\backslash\{0\}$ such that $\norm{x-y}\leq\eps\norm{X(x)}$ for $\eps<\eps_0$. For $t\in\frac{r_1}{2}\set{D}$, there exists a unique $u\in r_1\set{D}$ satisfying $g(t,u)=0$. Furthermore, $\Cmod{t-u}=O(\eps)$.

    Moreover, $\Jac_u(g)(t,u)\neq0$. In particular, by the implicit function theorem, there exists a $\class{1}$ function $f\colon\frac{r_1}{2}\set{D}\to r_1\set{D}$ such that $g(t,f(t))=0$ for $t\in\frac{r_1}{2}\set{D}$.
  \end{lem}

  \begin{pr} Let us note first that for $\Cmod{t}<r_0$, there exists a constant $C>1$ such that $C^{-1}\norm{x-y}\leq\norm{\flot{x}(t)-\flot{y}(t)}\leq C\norm{x-y}$. Indeed, by Gr\"{o}nwall Lemma, $C=e^{cr_0}$, for $c$ such that $\norm{X(x')-X(y')}\leq c\norm{x'-y'}$. Define $N(u)=\norm{\flot{x}(t)-\flot{y}(u)}^2$. By compactness of $\adhrD{r_1}$, $N$ attains its minimum somewhere in $u_0\in \adhrD{r_1}$. $N$ is of class $\class{1}$, so either $u_0$ satisfies $\bder{N}{u}(t_0)=g(t,u_0)=0$, or $\Cmod{u_0}=r_1$. We show that for $u$ too far from $t$, $N(u)>\norm{\flot{x}(t)-\flot{y}(t)}^2=N(t)$. Namely,
    \[\begin{aligned}\norm{\flot{x}(t)-\flot{y}(u)}&\geq\Cmod{t-u}\norm{X(\flot{y}(t))}-\norm{\flot{x}(t)-\flot{y}(t)}-\norm{\flot{y}(t)-\flot{y}(u)-(t-u)X(\flot{y}(t))},\\
        &\geq\left(\Cmod{t-u}+O(\eps)+O(\Cmod{t-u}^2)\right)\norm{X(\flot{y}(t))},\end{aligned}\]
    Note that $\norm{X(\flot{y}(t))}\geq C\norm{X(x)}$. It follows that for a sufficiently small $r_1$, $\eps_0$ sufficiently small chosen in consequence and $\Cmod{t-u}>O(\eps)$, $\norm{\flot{x}(t)-\flot{y}(u)}>C\eps\norm{X(x)}\geq\norm{\flot{x}(t)-\flot{y}(t)}$. Therefore, $t-u_0=O(\eps)$, $u_0$ does not belong to the boundary $\{\Cmod{u}=r_1\}$ and $g(t,u_0)=0$. We give in Figure \ref{figexistenceorthproj} a schematic view of the quantities involved in our computation.
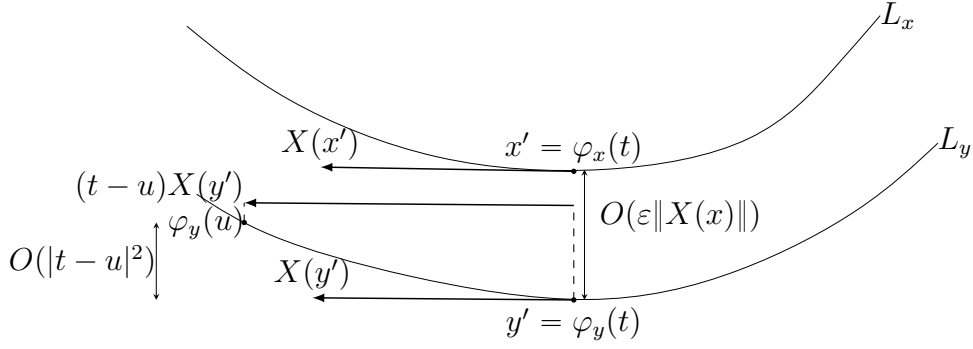
\begin{figure}
  \def \globalscale {10.000000}
\begin{tikzpicture}[y=0.80pt, x=0.80pt, yscale=-\globalscale, xscale=\globalscale, inner sep=0pt, outer sep=0pt]
\path[draw=black,line cap=butt,line join=miter,line width=0.056pt]
  (5.5579,8.0695) .. controls (7.0013,9.2648) and (8.4446,10.4600) ..
  (10.3806,11.5583) .. controls (12.3167,12.6565) and (14.7455,13.6577) ..
  (17.2917,14.2362) .. controls (19.8378,14.8148) and (22.5012,14.9706) ..
  (25.0545,14.8237) .. controls (27.6077,14.6767) and (30.0506,14.2269) ..
  (31.8553,13.3551) .. controls (33.6600,12.4833) and (34.8264,11.1896) ..
  (35.7414,10.1552) .. controls (36.6564,9.1208) and (37.3200,8.3457) ..
  (37.9836,7.5706) node[right] {$L_x$};

\path[draw=black,line cap=butt,line join=miter,line width=0.056pt]
  (6.0132,15.9930) .. controls (7.0222,16.6588) and (8.0311,17.3246) ..
  (9.6171,17.9767) .. controls (11.2032,18.6288) and (13.3663,19.2673) ..
  (15.9760,19.8536) .. controls (18.5857,20.4400) and (21.6418,20.9742) ..
  (24.3089,20.9582) .. controls (26.9760,20.9423) and (29.2538,20.3762) ..
  (31.3742,19.5704) .. controls (33.4946,18.7645) and (35.4573,17.7190) ..
  (36.9830,16.6806) .. controls (38.5086,15.6422) and (39.5971,14.6110) ..
  (40.6856,13.5798) node[right] {$L_y$};

\path[draw=black,line cap=butt,line join=miter,line width=0.56pt,->,>=latex]
  (23.6528,14.9126) -- (11.8002,14.7126) node[above,yshift=.1cm] {$X(x')$};

  \filldraw (23.6528,14.9126) node[above,yshift=.1cm] {$x'=\varphi_x(t)$} circle (.1);
  
\path[draw=black,line cap=butt,line join=miter,line width=0.56pt,->,>=latex]
 (23.6528,20.9790) -- (11.4217,20.8781) node[above,yshift=.1cm] {$X(y')$};

 \filldraw (23.6528,20.9790) node[below,yshift=-.1cm] {$y'=\varphi_y(t)$} circle (.1);

 \path[draw=black,line cap=butt,line join=miter,line width=0.56pt,<-,>=latex]
  (8.2343,16.3973) node[above left] {$(t-u)X(y')$} -- (23.6528,16.5241);

  \filldraw (8.2343,17.3273) node[left] {$\varphi_y(u)$} circle (.1);
  \draw[dashed] (23.6528,20.9790) -- (23.6528,16.5241);
  \draw[dashed] (8.2343,17.3273) -- (8.2343,16.3973);
\path[draw=black,line cap=butt,line join=miter,line width=0.056pt,<->,>=stealth]
(4.1237,17.3273) node[below left,yshift=-.3cm] {$O(\vert t-u\vert^2)$} -- (4.1237,20.9790);

\path[draw=black,line cap=butt,line join=miter,line width=0.056pt,<->,>=stealth]
(24.1528,14.8326) node[below right,yshift=-.4cm,xshift=.2cm] {$O(\varepsilon\Vert X(x)\Vert)$} -- (24.1528,20.9790);

\end{tikzpicture}
  \caption{Proof of Lemma \ref{existf}. Quantities involved in the proof that $u$ cannot be too far from $t$.\label{figexistenceorthproj}}
\end{figure}

    Now, we come to the uniqueness of $u_0$. Take another $u$ such that $g(t,u)=0$ and $\Cmod{u}<r_1$, and compute
    \[\begin{aligned} g(t,u)-g(t,u_0)=0=&\left\langle\flot{x}(t)-\flot{y}(u_0),X(\flot{y}(u))-X(\flot{y}(u_0))\right\rangle+\left\langle\flot{y}(u_0)-\flot{y}(u),X(\flot{y}(u))\right\rangle,\\
        0 =&\left((u_0-u)+O(\Cmod{u-u_0}^2)+O(\eps\Cmod{u-u_0})\right)\norm{X(\flot{y}(u))}^2.\end{aligned}\]
    For $r_1$ and $\eps$ small enough, this is possible only if $u=u_0$ and hence $u_0$ is unique. It remains to show that the Jacobian with respect to $u$ does not vanish in $(t,u_0)$. Derivating \eqref{defgxy} with respect to $u$ and $\bar{u}$, we get
    \begin{equation}\label{dergubaru} \der{g}{u}(t,u)=-\norm{X(\flot{y}(u))}^2,\quad\bder{g}{u}(t,u)=\left\langle\flot{x}(t)-\flot{y}(u),\dder{\flot{y}}{u}(u)\right\rangle.\end{equation}
    Using \eqref{bornederflot}, it follows that
    \begin{equation}\label{minJacug}\Jac_u(g)(t,u_0)=\begin{vmatrix} \disder{g}{u} & \cjg{\disbder{g}{u}} \\ \disbder{g}{u} & \cjg{\disder{g}{u}} \end{vmatrix}(t,u_0)\geq (1-C\eps)\norm{X(\flot{y}(u_0))}^4\geq C_2\norm{X(x)}^4>0.\end{equation}
    The existence and regularity of the function $f$ is a consequence of the implicit function theorem and the uniqueness of $u_0$ to glue together the local implicit functions.

    \qed
  \end{pr}

  We will need in Section~\ref{secmainthm} a $\class{2}-$estimate for the orthogonal projection. We begin by controlling the $\class{2}-$norm of $f$.

  \begin{lem}\label{norm2f} With the notations of Lemma \ref{existf}, the $\class{2}-$norm of $f$ satisfies
    \[\norm{f-\id}_{\class{2}}=O(\eps).\] 
  \end{lem}

  \begin{pr} Lemma \ref{existf} implies that $\norm{f-\id}_{\infty}=O(\eps)$. Let us now differentiate $f-\id$. The implicit function theorem relates the derivatives of $g$ and $f$ by
\begin{equation}\label{derf-id}\begin{aligned}\der{f}{t}(t)-1&=-\begin{vmatrix}\disder{g}{u} & \cjg{\disbder{g}{u}} \\ \disbder{g}{u} & \cjg{\disder{g}{u}} \end{vmatrix}^{-1}\begin{vmatrix}\disder{g}{t}+\disder{g}{u} & \cjg{\disbder{g}{u}} \\ \disbder{g}{u} & \cjg{\disder{g}{u}}\end{vmatrix}(t,f(t)),\\
     \bder{f}{t}(t)&=-\begin{vmatrix}\disder{g}{u} & \cjg{\disbder{g}{u}} \\ \disbder{g}{u} & \cjg{\disder{g}{u}} \end{vmatrix}^{-1}\left(\cjg{\der{g}{t}}\cdot\bder{g}{u}\right)(t,f(t)), \end{aligned}\end{equation}
 since $\bder{g}{t}=0$. By \eqref{dergubaru} and \eqref{bornederflot}, $\Cmod{\bder{g}{u}(t,f(t))}=O(\norm{\flot{x}(t)-\flot{y}(f(t))}\norm{X(x)})$. On the other hand,
 \begin{equation}\label{dergt+dergu} \der{g}{t}(t,u)+\der{g}{u}(t,u)=\left\langle X(\flot{x}(t))-X(\flot{y}(u)),X(\flot{y}(u))\right\rangle.\end{equation}
 It follows that $\Cmod{\der{g}{t}+\der{g}{u}}(t,f(t))$ is also $O(\norm{\flot{x}(t)-\flot{y}(f(t))}\norm{X(x)})$. Denote by $\delta_{xy}(t)=\flot{x}(t)-\flot{y}(f(t))$. By the previous estimates and \eqref{minJacug}, we have that $\norm{d(f-\id)}_{\infty}$ is a $O(\norm{\delta_{xy}}\norm{X(x)}^{-1})$. Since
  \[\norm{d_t\delta_{xy}}\leq\norm{d_t\left(\flot{x}-\flot{y}\right)}+\norm{d\flot{y}}_{\infty}\norm{d_t(f-\id)}\leq C\norm{\delta_{xy}(t)},\]
  Gr\"{o}nwall Lemma implies that
  %\begin{equation}\label{postGrondeltaxy} \norm{\delta_{xy}(t)}\leq\norm{x-\flot{y}(f(0))}e^{C\Cmod{t}}\leq C_3\norm{x-y},\quad t\in \frac{r_1}{2}\set{D}.\end{equation}
  \[\norm{\delta_{xy}(t)}\leq\norm{x-\flot{y}(f(0))}e^{C\Cmod{t}}\leq C_3\norm{x-y},\quad t\in \frac{r_1}{2}\set{D}.\]
  
  The estimate on the $\class{1}-$ norm follows. We continue with the derivatives of order 2 of $f-\id$. Recall that
  % \begin{equation}\label{difdet}d_A\det\cdot H=\Tr\left(\transp{\com(A)}H\right),\end{equation}
  \[d_A\det\cdot H=\Tr\left(\transp{\com(A)}H\right),\]
  where $\com(A)$ denotes the comatrix of $A$, $\Tr$ is the trace operator, and $d_A\det$ is the differential of the determinant with respect to the matrix $A$. Denote by $J(t)=\Jac_u(g)(t,f(t))$ and differentiate the holomorphic part of~\eqref{derf-id}.
  \[d_t\left(J\cdot\left(\der{f}{t}-1\right)\right)=-\Tr\left(\begin{pmatrix} \cjg{\disder{g}{u}} & -\cjg{\disbder{g}{u}} \\ -\disbder{g}{u} & \disder{g}{t}+\disder{g}{u}\end{pmatrix}d_{(t,f(t))}\begin{pmatrix}\disder{g}{t}+\disder{g}{u} & \cjg{\disbder{g}{u}} \\ \disbder{g}{u} & \cjg{\disder{g}{u}} \end{pmatrix}\right)\circ(dt,df).\]
  Here, we have denoted by $d_{(t,f(t))}\begin{pmatrix}\der{g}{t}+\der{g}{u} & \cjg{\bder{g}{u}} \\ \bder{g}{u} & \cjg{\der{g}{u}} \end{pmatrix}$ the differential of the obvious matrix application with respect to the point $(t,f(t))$. Note again that $\Cmod{\bder{g}{u}}$ and $\Cmod{\der{g}{t}+\der{g}{u}}$ are  $O\left(\norm{x-y}\norm{X(x)}\right)$, $\Cmod{\der{f}{t}-1}$ is a $O\left(\norm{x-y}{X(x)}^{-1}\right)$, and $\norm{g}_{\class{2}}\leq C\norm{X(x)}^2$ by \eqref{bornederflot}. Differentiating \eqref{dergt+dergu}, we get that
  \[\begin{aligned}\left(\der{}{t}+\der{}{u}\right)\left(\der{g}{t}+\der{g}{u}\right)&=\left\langle\dder{\flot{x}}{t}-\dder{\flot{y}}{u},X\circ\flot{y}\right\rangle,\\
      \left(\bder{}{t}+\bder{}{u}\right)\left(\der{g}{t}+\der{g}{u}\right)&=\left\langle X\circ\flot{x}-X\circ\flot{y},\dder{\flot{y}}{u}\right\rangle.\end{aligned}\]
  Hence, $\Cmod{\left(\der{}{t}+\der{}{u}\right)\left(\der{g}{t}+\der{g}{u}\right)}$ and also that $\Cmod{\left(\bder{}{t}+\bder{}{u}\right)\left(\der{g}{t}+\der{g}{u}\right)}$ are $O(\norm{x-y}\norm{X(x)})$. It follows that

  \[\norm{d\left(J\cdot\left(\der{f}{t}-1\right)\right)}_{\infty}\leq C\norm{X(x)}^2F(x,y),\]
  where $F$ is given by
  \[\begin{aligned}F(x,y)&=\Tr\left(\begin{pmatrix}\norm{X(x)} & \norm{x-y}\\\norm{x-y} & \norm{x-y}\end{pmatrix}\begin{pmatrix}\norm{x-y} & \norm{X(x)} \\ \norm{X(x)} & \norm{X(x)}\end{pmatrix}\right)+\norm{x-y}\norm{X(x)}\\
      &=5\norm{x-y}\norm{X(x)}.\end{aligned}\]
   By~\eqref{minJacug}, and using estimates on $\Cmod{\der{f}{t}-1}$ and $\norm{g}_{\class{2}}$, we get that $\norm{d\left(\der{f}{t}-1\right)}_{\infty}$ is a $O(\norm{x-y}\norm{X(x)}^{-1})$. Next, we differentiate the anti-holomorphic part of~\eqref{derf-id}.
  \[\bder{}{t}\left(J(t)\bder{f}{t}(t)\right)=\bder{g}{u}\cdot\bder{}{t}\left(\cjg{\der{g}{t}}(t,f(t))\right)+\cjg{\der{g}{t}}\left(\bder{f}{t}\cdot\dbder{g}{t}{u}+\cjg{\der{f}{t}}\cdot\bdder{g}{u}\right).\]

  Recall that $\norm{g}_{\class{2}}$ is a $O(\norm{X(x)}^2)$ by \eqref{bornederflot}, $\Cmod{\bder{f}{t}}$ is a $O(\norm{x-y}\norm{X(x)}^{-1})$, $\Cmod{\bder{g}{u}}(t,f(t))$ is a $O(\norm{x-y}\norm{X(x)})$ and $\Cmod{\bdder{g}{u}}(t,f(t))=\Cmod{\langle\flot{x}(t)-\flot{y}(f(t)),\frac{\partial^3\flot{y}}{\partial u^3}(f(t))\rangle}$ is a $O(\norm{x-y}\norm{X(x)})$. Similarly, the Jacobian on the left hand side does not change the estimate if it is inside or outside the derivative. Hence $\Cmod{\bdder{f}{t}}\leq C\norm{x-y}\norm{X(x)}^{-1}$. This and the same inequalities on $\norm{d\left(\der{f}{t}-1\right)}_{\infty}$ give us finally
  \[\norm{f-\id}_{\class{2}}\leq C_4\norm{x-y}\norm{X(x)}^{-1}=O(\eps).\]

  \qed

\end{pr}

\begin{figure}[hbt]
  \def \globalscale {10.000000}
\begin{tikzpicture}[y=0.80pt, x=0.80pt, yscale=-\globalscale, xscale=\globalscale, inner sep=0pt, outer sep=0pt]
\path[draw=black,line cap=butt,line join=miter,line width=0.65pt]
  (7.5491,27.2325) .. controls (7.6896,28.5035) and (8.2434,29.1837) ..
  (9.3015,29.8278) node[above right] {$L_x$} .. controls (10.5324,30.4362) and (12.6706,31.0971) ..
  (15.5774,31.4219) .. controls (18.513,31.8065) and (22.1436,31.7363) ..
  (24.8933,31.4989) .. controls (27.6429,31.2615) and (30.2815,30.8570) ..
  (33.1357,30.2002) .. controls (35.9899,29.5434) and (39.0596,28.6344) ..
  (42.1292,27.7253);

\path[draw=black,line cap=butt,line join=miter,line width=0.65pt]
  (8.7364,39.7193) .. controls (8.8412,38.5081) and (8.902,37.4066) ..
  (9.5355,36.7451) node[below right] {$L_y$}.. controls (10.1691,36.0835) and (11.5338,35.6963) ..
  (13.314,35.411) .. controls (15.2422,35.0875) and (17.6526,34.9006) ..
  (19.7579,34.7225) .. controls (21.8631,34.5444) and (24.2608,34.3942) ..
  (26.4894,34.3050) .. controls (28.7180,34.2157) and (30.7774,34.1874) ..
  (33.6111,33.5955) .. controls (36.4447,33.0036) and (40.0525,31.8482) ..
  (43.6604,30.6928);

\path[draw=black,line cap=butt,line join=miter,line width=0.156pt]
  (3.3685,31.2979) .. controls (3.3681,31.6876) and (3.3677,32.0773) ..
  (3.2370,32.4857) .. controls (3.1064,32.8940) and (2.8454,33.3211) ..
  (2.4429,33.6396) .. controls (2.0404,33.9580) and (1.4962,34.1678) ..
  (1.1748,34.2821) .. controls (0.8533,34.3964) and (0.7546,34.4152) ..
  (0.6558,34.4341);%haut-gauche 1

  \path[draw=black,line cap=butt,line join=miter,line width=0.156pt]
  (2.5685,30.4979) .. controls (2.5681,31.2876) and (2.5677,31.2773) ..
  (2.4370,31.6857) .. controls (2.3064,32.0940) and (2.0454,32.5211) ..
  (1.6429,32.8396) .. controls (1.2404,33.1580) and (0.6962,33.3678) ..
  (0.3748,33.4821) .. controls (0.0533,33.5964) and (-0.0454,33.6152) ..
  (-0.1442,33.6341);%haut-gauche 2

\path[draw=black,line cap=butt,line join=miter,line width=0.156pt]
(0.5401,35.3975) .. controls (0.7886,35.3317) and (1.2498,35.3619) ..
  (1.649,35.4091) .. controls (2.1447,35.4524) and (2.4638,35.6107) ..
  (2.8188,35.7993) .. controls (3.1738,36.033) and (3.5645,36.4873) ..
  (3.7334,36.8568) .. controls (3.9034,37.3291) and (3.9425,37.6193) ..
  (3.9869,38.0929);%bas-gauche 1

  \path[draw=black,line cap=butt,line join=miter,line width=0.156pt]
  (-0.4599,36.0975) .. controls (-0.2114,36.0317) and (0.2498,36.0619) ..
  (0.649,36.1091) .. controls (1.1447,36.1524) and (1.4638,36.3107) ..
  (1.8188,36.4993) .. controls (2.1738,36.733) and (2.5645,37.1873) ..
  (2.7334,37.5568) .. controls (2.9034,38.0291) and (2.9425,38.3193) ..
  (2.9869,38.7929);%bas-gauche2

\path[draw=black,line cap=butt,line join=miter,line width=0.156pt]
  (10.8370,26.4582) .. controls (11.9901,26.9568) and (13.1431,27.4555) ..
  (14.8368,27.8741) .. controls (16.5305,28.2928) and (18.7648,28.6315) ..
  (21.3205,28.6900) .. controls (23.8762,28.7486) and (26.7532,28.5270) ..
  (29.9939,27.8949) .. controls (33.2345,27.2627) and (36.8387,26.2201) ..
  (40.4429,25.1774);

\path[draw=black,line cap=butt,line join=miter,line width=0.156pt]
  (12.4119,40.5607) .. controls (14.0085,39.4476) and (15.6051,38.3344) ..
  (17.5774,37.6002) .. controls (19.5496,36.8660) and (21.8975,36.5106) ..
  (24.5885,36.3340) .. controls (27.2795,36.1573) and (30.3135,36.1593) ..
  (33.7125,35.6733) .. controls (37.1116,35.1873) and (40.8756,34.2134) ..
  (44.6396,33.2394);

\path[draw=black,line cap=butt,line join=miter,line width=0.156pt]
   (29.8836,30.9151) .. controls (30.0178,31.9797) and (30.1520,33.0442) ..
   (30.2862,34.1086);

%\path[draw=black,line cap=butt,line join=miter,line width=0.056pt]
 % (34.1068,30.0258) .. controls (34.4813,31.0893) and (34.8559,32.1528) ..
  %(35.2304,33.2162);

  \draw[line width=0.156pt] (34.1145,29.9834) --(35,33.25);

\path[draw=black,line cap=butt,line join=miter,line width=0.156pt]
  (30.2600,33.7843) -- (30.5601,33.7492) -- (30.5965,34.0579);

\path[draw=black,line cap=butt,line join=miter,line width=0.156pt]
  (34.9215,32.9604) -- (35.2112,32.8819) -- (35.2897,33.1715);

\path[draw=black,fill=black,fill
  opacity=1,line width=0.156pt] (31.8854,33.9) node[below,yshift=-0.1cm] {$y$} ellipse (0.003cm and
  0.003cm);

\path[draw=black,fill=black,fill
  opacity=1,line width=0.156pt] (29.8722,30.8770) node[above left,yshift=.1cm] {$x$} ellipse (0.003cm and
  0.003cm);

\path[draw=black,fill=black,fill
  opacity=1,line width=0.156pt] (30.2885,34.1052) node[below left,yshift=-.1cm] {$\Phi_{xy}(x)$} ellipse (0.003cm and
  0.003cm);

\path[draw=black,fill=black,fill
  opacity=1,line width=0.156pt] (35,33.26) node[below right] {$\varphi_y(f(t))=\Phi_{xy}(x')$} ellipse (0.003cm and
  0.003cm);

\path[draw=black,fill=black,fill
  opacity=1,line width=0.156pt] (34.1145,29.9834) node[above,yshift=.1cm] {$x'=\varphi_x(t)$} ellipse (0.0030cm and
  0.0030cm);

\path[draw=black,fill=black,line width=0.156pt,miter limit=4.00]
  (3.8833,34.4790) node[right,xshift=-0.3cm] {$0$} ellipse (0.0060cm and 0.0060cm);

  \path[draw=black,line cap=butt,line join=miter,line width=0.156pt]
  (4.7834,38.5547) .. controls (4.8712,37.6051) and (4.9844,36.9194) ..
  (5.1187,36.378) .. controls (5.2281,35.8838) and (5.3483,35.5181) ..
  (5.8124,35.2128) .. controls (6.2763,34.9075) and (7.2165,34.5781) ..
  (7.976,34.3722) .. controls (8.8028,34.1398) and (9.3779,34.0739) ..
  (10.3060,33.9779);%bas-droite1

  \path[draw=black,line cap=butt,line join=miter,line width=0.156pt]
  (6.1347,38.6678) .. controls (6.2771,37.4938) and (6.2108,36.5668) ..
  (6.6271,36.0372) .. controls (7.0435,35.5075) and (8.1594,35.2721) ..
  (9.4321,34.956);%bas-droite 2
  %.. controls (11.82095,35.09895) and (13.2401,34.75585) ..
  %(14.62855,34.534);
  %.. controls (16.017,34.31215) and (17.3746,34.21155) ..
  % (18.6477,34.14145);

\path[draw=black,line cap=butt,line join=miter,line width=0.156pt]
  (4.3016,30.6884) .. controls (4.5509,31.5799) and (4.7532,32.3575) ..
  (5.03364,32.7973) .. controls (5.2745,33.1961) and (5.6703,33.371) ..
  (6.3741,33.3722) .. controls (7.1547,33.384) and (8.4376,33.2151) ..
  (9.9145,33.0391);% Haut-droite 1

    \path[draw=black,line cap=butt,line join=miter,line width=0.156pt]
  (5.0472,29.7325) .. controls (5.4966,30.6216) and (5.7435,31.6504) ..
  (6.4879,32.1578) .. controls (7.2323,32.6651) and (8.8463,32.4545) ..
  (10.1827,32.4731);%haut-droite 2

\end{tikzpicture}
  \caption{Orthogonal projection from $\leafu{x}$ near $x$ onto $\leafu{y}$ near $y$ (outside the singular point $0$)\label{picorthproj}}
\end{figure}
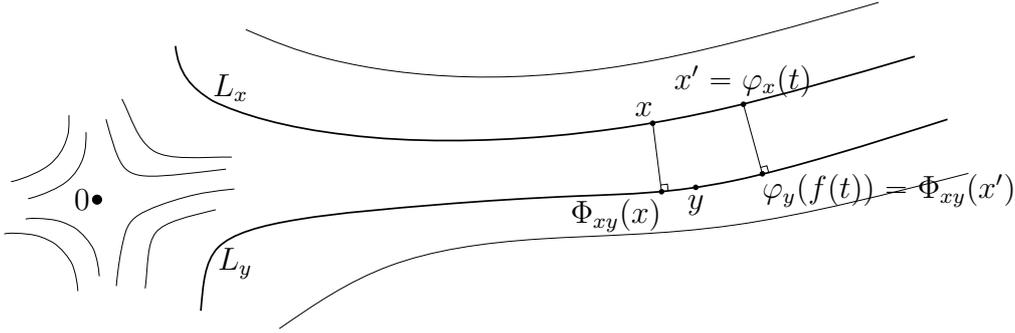

Define then locally $\Phi_{xy}=\flot{y}\circ f\circ\flot{x}^{-1}\colon\germ{\leafu{x}}{x}\to\germ{\leafu{y}}{y}$. By definition of $f$, $\Phi_{xy}$ is the local orthogonal projection from $\leafu{x}$ in the neighbourhood of $x$ to $\leafu{y}$ in the neighbourhood of $y$. Figure \ref{picorthproj} gives a schematic view of this definition and its link with the function $f$. The radii of these neighbourhoods are $O\left(\norm{x}\right)$ by \eqref{borneX} and Lemma~\ref{existf}. Considering the local definition of $\Phi_{xy}$, \eqref{borneX} and Lemma \ref{norm2f}, we get easily that
  \begin{equation}\label{norm2Phi}\norm{\Phi_{xy}-\id}_{\infty}\leq C_5\norm{x-y},\quad\norm{\Phi_{xy}-\id}_{\class{1}}\leq C_5\frac{\norm{x-y}}{\norm{x}},\quad\norm{\Phi_{xy}-\id}_{\class{2}}\leq C_5\frac{\norm{x-y}}{\norm{x}^2}.\end{equation}

  We will also need some estimates on the distances induced by the standard Hermitian metric restricted to $\leafu{x}$ and $\leafu{y}$. Denote by $\dhermf{}$ the distance induced by the restriction of the standard Hermitian metric on a leaf $\leaf{}$. Let $x_1$ and $x_2$ be two points in a sufficiently small neighbourhood $U$ of $x$ in $\leafu{x}$ such that $\flot{x}$ is bijective from a neighbourhood of $0$ onto $2U$. This neighbourhood can be chosen of radius $O(\norm{x})$ by \eqref{borneX} and \eqref{bornederflot}. Suppose that $x_1=\flot{x}(t_1)$ and $x_2=\flot{x}(t_2)$ and denote by $\gamma$ a path of minimal length in $\leafu{x}$ from $x_1$ to $x_2$. The path $\gamma$ stays in $2U$. Then, $\flot{y}\circ f\circ\flot{x}^{-1}\circ\gamma$ is a path in $\leafu{y}$ from $\Phi_{xy}(x_1)$ to $\Phi_{xy}(x_2)$. It follows easily that there exists a constant $C_6$ such that
  \begin{equation}\label{bounddgLyPhixy} \dhermf{y}(\Phi_{xy}(x_1),\Phi_{xy}(x_2))\leq C_6\dhermf{x}(x_1,x_2).\end{equation}

  \section{Proof of the main theorem \label{secmainthm}}

  We are now ready to prove Theorem~\ref{mainthm}. Throughout this section, we will have to consider three different distances and will be very careful to distinguish them by our notations. We denote by $\dherm{}$ the distance induced by $g_{\mani{M}}$ on $\mani{M}$. For $\leaf{}$ a leaf of $\fol{}$, we consider two distances:
  \begin{itemize}
  \item $\dPC{}$ the distance induced by the Poincar\'e metric on $\leaf{}$ ;
  \item $\dhermf{}$, as in Section~\ref{secorthproj}, the distance induced by the restriction of $g_{\mani{M}}$ on $\leaf{}$.
  \end{itemize}
  For a set $K$ and two functions $u,v$ from $K$ to $\leaf{}$ or $\mani{M}$, we define the distances between $u$ and $v$.

\[\dhermset{K}(u,v)=\underset{x\in K}{\sup}\,\dherm{}(u(x),v(x)),\quad\dPCset{K}(u,v)=\underset{x\in K}{\sup}\,\dPC{}(u(x),v(x)),\]
\[\dhermfset{}{K}(u,v)=\underset{x\in K}{\sup}\,\dhermf{}(u(x),v(x)).\]

By compactness of $M$, it is sufficient to show the result~\eqref{eqmainthm} only for $x$ and $y$ close to each other. Hence, we will suppose that
\[\frac{\max\left(\log^{\star}\dsing{E}{x},\log^{\star}\dsing{E}{y}\right)}{\log^{\star}d(x,y)}\leq e^{-\alpha^{-1}R},\]
for $R$ big enough and some $\alpha$ that we will specify later.

  \subsection{Geometric setup\label{subsecgeomsetup}}

   Note that the existence of $C$ in~\eqref{eqmainthm} is independent of the metric $g_{\mani{M}}$. Then, we may build $g_{\mani{M}}$ by a partition of the unity satisfying some peculiar hypothesis that will help us in our proof. By compactness of $\mani{M}$, we can find a finite open covering $~\mathcal{U}=\left(U_p,U_i\right)_{p\in E,i\in I}$, where
  \begin{enumerate}
  \item $U_p$ is a neighbourhood of $p$ and we have a holomorphic chart $U_p\simeq4\rho\set{D}^n$ centered at $p$ in which $g_{\mani{M}}$ is the standard Hermitain metric $\norm{dz}^2$ on a smaller polydisk $3\rho\set{D}^n$. We suppose that $\fol{}$ is generated on $U_p$ by some holomorphic vector field $X_p$.
  \item $\left(U_i\right)_{i\in I}$ is an open covering of $\mani{M}\backslash\bigcup_{p\in\mani{E}}\frac{1}{4}U_p$ by flow boxes. We also suppose that $\fol{}$ is generated on $U_i$ by a holomorphic vector field $X_i$.
  \end{enumerate}

  To simplify the computations, we suppose that the diameter of $\mani{M}$ for $g_{\mani{M}}$ is lower than 1. That way, for $\delta$ a distance on $\mani{M}$, $\log^{\star}\delta=1-\log\delta$. Denote by $d>0$ some small distance such that any point of $\mani{M}$ admits a holomorphic neighbourhood of radius at least $d$. We suppose that $\fol{}$ is Brody hyperbolic and we will denote by $A$ the constant appearing in Definition~\ref{defBrody}.

  The following important local estimate on $\eta$ is needed.

  \begin{thm}[Lins Neto-Canille Martins~{\cite[Theorem~3]{MLN}}]\label{thmMLN} Let $\fol$ be a holomorphic foliation by Riemann surfaces on a ball $r\set{B}$ with non-degenerate singularity at the origin, where $\set{B}$ is the unit ball of $\set{C}^n$. Then there exists a smaller radius $0<\rho<r$ such that
    \[C^{-1}\norm{z}\log^{\star}\norm{z}<\eta(z)<C\norm{z}\log^{\star}\norm{z},\quad z\in \rho\set{B}\backslash\{0\}.\]
  \end{thm}

  We obtain the following global estimate on $\eta$.
  
    \begin{prop} \label{setupbndlamb} There exists a constant $C_1>1$ such that
    \[C_1^{-1}\dsing{E}{x}\log^{\star}\dsing{E}{x}\leq\eta(x)\leq C_1\dsing{E}{x}\log^{\star}\dsing{E}{x},\quad x\in\manis{M}{E}.\]
    \end{prop}

    \begin{pr} Suppose first that $x\in\manis{M}{E}$ is far from $E$. By the embedding of the plaque passing through $x$ in a flow box, we get $\eta(x)>C$. Moreover, by Brody hyperbolicity of $\fol{}$, we have $\eta(x)<A$. Since $x$ is far from $E$, the conclusion of Proposition~\ref{setupbndlamb} follows.

      Suppose next that $x\in \frac{3}{4}U_p$ for some $p\in E$. Considering the restriction of $\fol{}$ to $U_p$ and applying Theorem~\ref{thmMLN}, we get a holomorphic function $u\colon\set{D}\to\leafu{x}\cap U_p$ such that $u(0)=x$ and $\eta(x)\geq\norm{u'(0)}_{g_{\mani{M}}}>C^{-1}\dsing{E}{x}\log^{\star}\dsing{E}{x}$. For the upper estimate, recall that $\fol{}$ is Brody hyperbolic. Hence, there exists a radius $R_0$, independant of $x\in\frac{3}{4}U_p$, such that the disk in $\leafu{x}$ of radius $R_0$ with respect to the Poincar\'{e} metric is contained in $U_p$. Denote by $u_x$ a uniformization of $\leafu{x}$ such that $u_x(0)=x$. For the corresponding Euclidean radius $r_0=\frac{e^{R_0}-1}{e^{R_0}+1}$, we have $u_x(\rD{r_0})\subset U_p$. By the extremal property of Poincar\'{e} metric and Theorem~\ref{thmMLN}, we get $\eta(x)<r_0^{-1}C\dsing{E}{x}\log^{\star}\dsing{E}{x}$.
      \qed
      \end{pr}
      
      We are going to prove our H\"{o}lder estimate by building a nearly holomorphic mapping $\psi\colon\DR{R}\to\leafu{y}$, which is close to a uniformization of $\leafu{x}$. Next, we solve a Beltrami equation to have a holomorphic function close to $\psi$ which will make $\eta(x)$ and $\eta(y)$ close. To produce this nearly holomorphic mapping, we shall use the orhogonal projection we studied earlier. Note that in~\cite[Lemmas~2.11 and~2.12]{DNSII}, the three authors use the invariance by homothety of $\fol{}$ in the case of linearizable singularities. Here, we had to establish in Section~\ref{secorthproj} the same estimates for orthogonal projections, since this trick is not available anymore in our context.

      Denote by $\ballleafu{x}{r}$ the ball of center $x$ and radius $r$ for the distance $\dhermf{x}$. What we have proven in Section~\ref{secorthproj} (see \eqref{norm2Phi}, \eqref{bounddgLyPhixy} and the fact that $\Phi_{xy}$ is defined on a neighourhood of radius $O(\norm{x})$), together with the classical estimates of orthogonal projection on non-singular flow boxes, can be summarized by the following lemma in our geometric setup.

      \begin{lem}\label{geomsetuporthproj} There exist constants $\eps_0$, $\eps_1$, $k$ and $K$ such that for $x,y\in\manis{M}{E}$, if $\dherm{}(x,y)\leq\eps_1\dsing{E}{x}$, then there exists a local orthogonal projection
        \[\Phi_{xy}\colon\ballleafu{x}{\eps_0\dsing{E}{x}}\to\ballleafu{y}{k\eps_0\dsing{E}{y}},\]
        satisfying
        \begin{enumerate}
        \item $\dhermf{y}(y,\Phi_{xy}(x))\leq k\dherm{}(x,y)$,
        \item for $x_1,x_2\in\ballleafu{x}{\eps_0\dsing{E}{x}}$, $\dhermf{y}(\Phi_{xy}(x_1),\Phi_{xy}(x_2))\leq k\dhermf{x}(x_1,x_2).$
        \item $\Phi_{xy}$ is smooth and we have the following estimates for $\Phi_{xy}-\id$ and its derivatives in the charts $U_{i}$ and $U_p$.
        \[\norm{\Phi_{xy}-\id}_{\infty}\leq e^K\dherm{}(x,\Phi_{xy}(x)),\quad\norm{\Phi_{xy}-\id}_{\class{1}}\leq e^K\frac{\dherm{}(x,\Phi_{xy}(x))}{\dsing{E}{x}},\]
        \[\norm{\Phi_{xy}-\id}_{\class{2}}\leq e^K\frac{\dherm{}(x,\Phi_{xy}(x))}{\dsing{E}{x}^2}.\]
      \item If $x'\in\ballleafu{x}{\eps_0\dsing{E}{x}}$, $y'\in\ballleafu{y}{k\eps_0\dsing{E}{x}}$ and $\dherm{}(x',y')\leq\eps_1\dsing{E}{x'}$, then $\Phi_{x'y'}=\Phi_{xy}$ on the intersection of their domains of definition.
      \end{enumerate}
      \end{lem}
      
  \subsection{Proof of the Hölder estimate}

  We begin our proof with the second step. Take $x,y\in\manis{M}{E}$ and take $u_x\colon\set{D}\to\leafu{x}$ a uniformization of $\leafu{x}$ such that $u_x(0)=x$. Such a covering is unique up to a rotation. Take also $u_y\colon\set{D}\to\leafu{y}$ a uniformization of $\leafu{y}$ such that $u_y(0)=y$. The following lemma is very similar to~\cite[Proposition~3.6]{DNSII}, the proof of which is partially referred to~\cite[Proposition~2.2]{DNSI}.
  
  \begin{lem}\label{Rdelimpsemilamb} Let $R>0$ and suppose that
    \begin{enumerate}
    \item \label{lemRdelpsi} There exists a smooth map $\psi\colon\adhDR{R}\to\leafu{y}$ without critical point such that $\psi(0)=y$ and $\dhermset{\adhDR{R}}\left(\psi,u_x\right)\leq e^{-2R}$. \label{existspsi}
    \item \label{lemRdeldpsi} $\norm{d\psi}_{\infty}\leq 2A$ for the constant $A$ appearing in Definition~\ref{defBrody}. We consider the norm of $d\psi$ with respect to the Poincar\'e metric on the source $\adhDR{R}$ and the restriction of $g_{\mani{M}}$ on the goal $\leafu{y}$.
    \item \label{lemRdelPsi} There exists a smooth map $\Psi\colon\adhDR{R}\to\set{D}$ such that $\Psi(0)=0$, $u_y\circ\Psi=\psi$ and the Beltrami coefficient $\mu_{\Psi}$ of $\Psi$ satisfies $\norm{\mu_{\Psi}}_{\class{1}}\leq e^{-2R}$. Here, we consider the norm of $\mu_{\Psi}$ with respect to the standard Euclidean metric on $\set{D}$ at the goal and the source.
    \end{enumerate}

    If $R$ is sufficiently large, then there exists a constant $C_2$ such that
      \[\eta(x)-\eta(y)\leq C_2e^{-R}.\]
    \end{lem}

We recall that the Beltrami coefficient $\mu_{\Psi}$ is defined by $\bder{\Psi}{t}=\mu_{\Psi}\der{\Psi}{t}$.

\begin{rem} Note that the existence of $\Psi$ in~\eqref{lemRdelPsi} such that $\Psi(0)=0$ and $u_y\circ\Psi=\psi$ is a consequence of~\eqref{lemRdelpsi} and the fact that $\adhDR{R}$ is simply connected. The smooth map $\Psi$ is just the unique lifting of $\psi$ \emph{via} $u_y$ such that $\Psi(0)=0$.
\end{rem}

    \begin{pr} Let us first translate our conditions in a commutative diagram. The mention $\hol{}$ (resp. $\smooth$) on an arrow will denote a holomorphic (resp. smooth) function.
      \[\begin{gathered}\xymatrixcolsep{4pc}\xymatrix{& \set{D}\ar[d]^{u_y}_{\hol}\\\adhDR{R}\ar[ur]^{\Psi}_{\smooth}\ar[r]^{\smooth}_{\psi}\ar@{^{(}->}[d]& \leafu{y}\\\set{D}\ar[r]^{\hol}_{u_x}&\leafu{x}}\end{gathered}.\]

      We want to find $v\colon\DR{R}\to\set{D}$ holomorphic and $q\colon\adhDR{R}\to\adhDR{R}$ smooth and close to the identity such that $\Psi=v\circ q$. To make $v$ holomorphic, we should have $\bder{q}{t}=\mu_{\Psi}\der{q}{t}$. It was proven by Schatz and Earle in~\cite{Teich} that this Beltrami equation can be solved with $q(0)=0$, and $\norm{q-\id}_{\infty}\leq \kappa\norm{\mu_{\Psi}}_{\class{1}}\leq\kappa e^{-2R}$, for some constant $\kappa>0$. We also have $\norm{q^{-1}-\id}_{\infty}\leq\kappa e^{-2R}$ for another constant $\kappa$. By construction, $v=\Psi\circ q^{-1}$ is holomorphic. Let us complete our commutative diagram.
      %\begin{equation}\label{comdiagRdelsemi}\begin{gathered}\xymatrixcolsep{4pc}\xymatrix{\adhDR{R}\ar[r]^v_{\hol}& \set{D}\ar[d]^{u_y}_{\hol}\\\adhDR{R}\ar[ur]^{\Psi}_{\smooth}\ar[r]^{\smooth}_{\psi}\ar@{^{(}->}[d]\ar[u]^q_{\smooth}& \leafu{y}\\\set{D}\ar[r]^{\hol}_{u_x}&\leafu{x}}\end{gathered}\end{equation}
      \[\begin{gathered}\xymatrixcolsep{4pc}\xymatrix{\adhDR{R}\ar[r]^v_{\hol}& \set{D}\ar[d]^{u_y}_{\hol}\\\adhDR{R}\ar[ur]^{\Psi}_{\smooth}\ar[r]^{\smooth}_{\psi}\ar@{^{(}->}[d]\ar[u]^q_{\smooth}& \leafu{y}\\\set{D}\ar[r]^{\hol}_{u_x}&\leafu{x}}\end{gathered}\]

      The estimate for $\eta$ will come from a comparison between the two holomorphic mappings $u_y\circ v$ and $u_x$. Denote by $v_y\colon\set{D}\to\leafu{y}$ defined by $v_y(t)=u_y\circ v(rt)$ with $r=\frac{e^R-1}{e^R+1}$ the Euclidean radius associated to $R$. Since $v_y$ is holomorphic and $v_y(0)=\psi\circ q^{-1}(0)=y$, by definition of $\eta$, we have $\eta(y)\geq\norm{v_y'(0)}_{g_{\mani{M}}}$. Consider a radius $r_0$ such that $u_x\left(\adhrD{r_0}\right)$ and $v_y\left(\adhrD{r_0}\right)$ are both included in a common holomorphic chart of $\mani{M}$. It is clear that this radius can be chosen at least equal to $\frac{d-e^{-2R}}{2A}\geq\frac{d}{4A}$ if $R$ is large enough. We apply the Cauchy formula in this chart. We get
      \begin{equation}\label{Cauchycomplambda}\norm{u_x'(0)-v_y'(0)}_{g_{\mani{M},x}}\leq C\norm{u_x'(0)-v_y'(0)}\leq\frac{4AC}{d}\underset{t\in\rD{r_0}}{\sup}\norm{u_x(t)-v_y(t)}\leq\frac{4CC'A}{d}\underset{t\in\rD{r_0}}{\sup}\dherm{}(u_x(t),v_y(t)),\end{equation}
      where the norm without index is the Euclidean norm in the chart, and the constants $C, C'$ are given by the equivalence of Hermitian metrics. Note also that we have
      \begin{equation}\label{complambda}\eta(x)-\eta(y)\leq\norm{u_x'(0)}_{g_{\mani{M},x}}-\norm{v_y'(0)}_{g_{\mani{M},y}}\leq\norm{u_x'(0)-v_y'(0)}_{g_{\mani{M},x}}+C''A\dherm{}(x,y),\end{equation}
      where the constant $C''$ comes from the difference between $\norm{v_y'(0)}_{g_{\mani{M},x}}$ and $\norm{v_y'(0)}_{g_{\mani{M},y}}$, that is, from the spatial variation of $g_{\mani{M}}$. Since $\dherm{}(x,y)\leq e^{-2R}$, the second term is well controlled.
      
      By~\eqref{complambda} and~\eqref{Cauchycomplambda}, it is sufficient to estimate $\dhermset{\rD{r_0}}(u_x,v_y)$. Take $t\in\rD{r_0}$ and remember that $u_y\circ v=\psi\circ q^{-1}$. We have
      \begin{equation}\label{estimateCauchycomplambda}\begin{aligned}\dherm{}(u_x(t),v_y(t))&\leq\dherm{}(u_x(t),u_x(rt))+\dherm{}(u_x(rt),\psi(rt))+\dherm{}\left(\psi(rt),\psi\left(q^{-1}(rt)\right)\right)\\
          &\leq(1-r)r_0\underset{\tau\in\rD{r_0}}{\sup}\norm{u_x'(\tau)}_{g_{\mani{M}}}+e^{-2R}+\norm{d\psi}_{\infty}\dPC{}\left(rt,q^{-1}(rt)\right).\end{aligned}\end{equation}

      Using an automorphism of $\set{D}$, we show that for $\Cmod{\tau}\leq r_0$, $\norm{u_x'(\tau)}_{g_{\mani{M}}}=\frac{\norm{u_{u_x(\tau)}'(0)}_{g_{\mani{M}}}}{1-\Cmod{\tau}^2}\leq\frac{A}{1-r_0^2}$. Moreover, $1-r=\frac{2}{e^R+1}\leq2e^{-R}$. Since $\dherm{}\left(q^{-1}(rt),rt\right)\leq Ce^{-2R}$, by Lemma~\ref{lemPCeucl}, $\dPC{}\left(rt,q^{-1}(rt)\right)\leq Ce^{-R}$. Back to~\eqref{estimateCauchycomplambda}, we get some constant $C'$ such that
      \[\dherm{}(u_x(t),v_y(t))\leq C'e^{-R}.\]

      Finally, by~\eqref{Cauchycomplambda} and~\eqref{complambda}, we obtain
      \[\eta(x)-\eta(y)\leq C_2e^{-R}.\]
      \qed
    \end{pr}

    We continue the proof of Theorem~\ref{mainthm} similarly to~\cite[Proposition~3.7]{DNSII}. We show that if $x$ and $y$ are sufficiently close, then there exist $\tilde{y}\in\leafu{y}$ close to $y$ and $\tilde{x}\in\leafu{x}$ close to $x$ such that $x,\tilde{y}$ and $y,\tilde{x}$ satisfy the hypothesis of Lemma~\ref{Rdelimpsemilamb}. We suppose more precisely that
    \begin{equation}\label{hypmainthm}\log^{\star}\dherm{}(x,y)\geq e^{\alpha^{-1}R}\max\left(\log^{\star}\dsing{E}{x},\log^{\star}\dsing{E}{y}\right),\end{equation}
    for some fixed $\alpha$ we will specify later. We begin by some work on $\leafu{x}$ which will be helpful to build the function $\psi$ by orthogonal projection.

    Consider $\xi\in\set{D}$ such that $\dPC{}(0,\xi)=R$. We want to subdivide the geodesic $\left[0,\xi\right]$ in subsegments $\left[\xi_j,\xi_{j+1}\right]$, for $j\in\llbracket0,N-1\rrbracket$. Set $\xi_0=0$ and for $j\in\llbracket0,N\rrbracket$, $x_j=u_x(\xi_j)$. Choose at each step $\xi_{j+1}\in\inter{\xi_j}{\xi}$ such that for all $\zeta\in\left[\xi_j,\xi_{j+1}\right]$, $\dhermf{x}\left(u_x(\zeta),x_j\right)\leq\eps_0\dsing{E}{x_j}$ and $\dhermf{x}(x_{j+1},x_j)=\eps_0\dsing{E}{x_j}$. At last step, we set $\xi_N=\xi$ and suppose that for all $\zeta\in\left[\xi_{N-1},\xi_N\right]$, $\dhermf{x}\left(u_x(\zeta),x_{N-1}\right)\leq\eps_0\dsing{E}{x_{N-1}}$. Denote by $R_j=\dPC{}(0,\xi_j)$ and $r_j=\Cmod{\xi_j}$. This setup will enable us to consider the local orthogonal projections $\Phi_{x_jy_j}$ by Lemma \ref{geomsetuporthproj}, where $y_0=y$ and $y_{j+1}=\Phi_{x_jy_j}(x_{j+1})$ (see Lemma \ref{lemconstructionpsi} below).

    At this stage, it is not certain that this process stops and that $N$ is finite, but the lemmas below will give us a bound for $N$.

    \begin{lem}\label{bndRj+1-Rj} We have the following estimate comparing $\dhermf{x}(x_j,x_{j+1})$ and $\dPC{}(\xi_j,\xi_{j+1})$.
      \[\dhermf{x}(x_j,x_{j+1})\leq\frac{1}{2}\underset{\left[\xi_j,\xi_{j+1}\right]}{\sup}\left(\eta\circ u_x\right)\,\dPC{}(\xi_j,\xi_{j+1})=\frac{1}{2}\underset{\left[\xi_j,\xi_{j+1}\right]}{\sup}\left(\eta\circ u_x\right)\,(R_{j+1}-R_j).\]
    \end{lem}

    \begin{pr} Denote by $\gamma$ a parametrization of the geodesic $\left[\xi_j,\xi_{j+1}\right]$. We have
      \[\begin{aligned} \dhermf{x}(x_j,x_{j+1})&\leq\int_0^1\norm{\left(u_x\circ\gamma\right)'(t)}_{g_{\mani{M}}}dt\\
          &\leq\underset{\left[\xi_j,\xi_{j+1}\right]}{\sup}\left(\eta\circ u_x\right)\,\int_0^1\frac{\Cmod{\gamma'(t)}\norm{u_x'(\gamma(t))}_{g_{\mani{M}}}}{\eta(u_x(\gamma(t)))}dt.\end{aligned}\]
      Now, $\norm{u_x'(\gamma(t))}_{g_{\mani{M}}}=\frac{\norm{u_{u_x(\gamma(t))}'(0)}_{g_{\mani{M}}}}{1-\Cmod{\gamma(t)}^2}=\frac{\eta(u_x(\gamma(t)))}{1-\Cmod{\gamma(t)}^2}$.
      It follows that
      \[\begin{aligned}\dhermf{x}(x_j,x_{j+1})&\leq\underset{\left[\xi_j,\xi_{j+1}\right]}{\sup}\left(\eta\circ u_x\right)\,\int_0^1\frac{\Cmod{\gamma'(t)}}{1-\Cmod{\gamma(t)}^2}dt\\
          &\leq\frac{1}{2}\underset{\left[\xi_j,\xi_{j+1}\right]}{\sup}\left(\eta\circ u_x\right)\,\dPC{}(\xi_j,\xi_{j+1}).\end{aligned}\]

      Since $\dPC{}(\xi_j,\xi_{j+1})=R_{j+1}-R_j$, the result of the lemma follows.
      \qed
    \end{pr}

    The following lemma is an analog of~\cite[Lemma~3.4]{DNSII} in our context.
    
    \begin{lem}\label{nottoocloseS} Let $\zeta,\sigma$ be points of $\set{D}$. Denote by $z=u_x(\zeta)$, $w=u_x(\sigma)$. There exist constants $C_3, C_4\in\set{R}_+^*$, such that we have
      \[\Cmod{\ln\frac{\log^{\star}\dsing{E}{z}}{\log^{\star}\dsing{E}{w}}}\leq C_3+C_4\dPC{}(\zeta,\sigma).\]
    \end{lem}

    \begin{pr} Since $\zeta$ and $\sigma$ play a symmetric role, assume that $\dsing{E}{z}\geq\dsing{E}{w}$. We have
      \begin{equation}\label{Cmodlnfracloglog}\begin{aligned}\Cmod{\ln\frac{\log^{\star}\dsing{E}{z}}{\log^{\star}\dsing{E}{w}}}&=\ln\frac{\log^{\star}\dsing{E}{w}}{\log^{\star}\dsing{E}{z}}=-\ln\left(\log^{\star}\dsing{E}{z}\right)+\ln\left(\log^{\star}\dsing{E}{w}\right)\\
          &=C\int_{\dsing{E}{w}}^{\dsing{E}{z}}\frac{dt}{t\log^{\star}t}.\end{aligned}\end{equation}

      Let $\gamma\colon\inter{0}{1}\to\set{D}$ be a parametrization of the geodesic from $\sigma$ to $\zeta$. Consider $\tau=\sup\left\{t\in\inter{0}{1}\setst{}\forall t'\in\inter{0}{t},\,\dsing{E}{u_x(\gamma(t'))}\leq\min\left(\dsing{E}{z},3\rho\right)\right\}$, with by convention $\tau=0$ if $\dsing{E}{w}>3\rho$. Denote by $\zeta'=\gamma(\tau)$ and $w'=u_x(\zeta')$. For $t\in\inter{0}{\tau}$, $u_x(\gamma(t))$ is contained in a singular chart $\frac{3}{4}U_p$ on which $g_{\mani{M}}$ is the standard Hermitian metric of $\set{C}^n$. Therefore $\dsing{E}{u_x(\gamma(t))}=\norm{u_x(\gamma(t))}$ and the second inequality of Proposition~\ref{setupbndlamb} will apply on it. We split our integral in two.
      \[\int_{\dsing{E}{w}}^{\dsing{E}{z}}\frac{dt}{t\log^{\star}t}=\int_{\dsing{E}{w}}^{\dsing{E}{w'}}\frac{dt}{t\log^{\star}t}+\int_{\dsing{E}{w'}}^{\dsing{E}{z}}\frac{dt}{t\log^{\star}t}.\]

      Now, the second integral is either on the empty set if $\dsing{E}{z}\leq3\rho$, or we have $\dsing{E}{w'}=3\rho$ and $\dsing{E}{z}\leq 1$. In both cases, the second integral is bounded above by $\int_{3\rho}^1\frac{dt}{t\log^{\star}t}$. We treat the first integral by using the standard Hermitian metric on $\set{C}^n$ and a change of variable. In the following, $\Re z$ denotes the real part of a complex number $z$.
      \[\begin{aligned}\int_{\dsing{E}{w}}^{\dsing{E}{w'}}\frac{dt}{t\log^{\star}t}&=\int_{\norm{u_x(\gamma(0))}}^{\norm{u_x(\gamma(\tau))}}\frac{dt}{t\log^{\star}t}=\int_0^{\tau}\frac{\Re \left\langle\left(u_x\circ\gamma\right)'(t),u_x\circ\gamma(t)\right\rangle}{\norm{u_x(\gamma(t))}^2\log^{\star}\norm{u_x(\gamma(t))}}dt\\
            &\leq\int_0^{\tau}\frac{\norm{u_x'(\gamma(t))}\Cmod{\gamma'(t)}}{\norm{u_x(\gamma(t))}\log^{\star}\norm{u_x(\gamma(t))}}dt\\
            &\leq C_1\int_0^{\tau}\frac{\norm{u_x'(\gamma(t))}\Cmod{\gamma'(t)}}{\eta(u_x(\gamma(t)))}dt=C_1\int_0^{\tau}\frac{\Cmod{\gamma'(t)}}{1-\Cmod{\gamma(t)}^2}dt\\
            &\leq\frac{C_1}{2}\dPC{}(\zeta,\zeta')\leq\frac{C_1}{2}\dPC{}(\zeta,\sigma).\end{aligned}\]

        Coming back to~\eqref{Cmodlnfracloglog}, and using our work on the second integral, we obtain
        \[\Cmod{\ln\frac{\log^{\star}\dsing{E}{z}}{\log^{\star}\dsing{E}{w}}}\leq C\int_{3\rho}^1\frac{dt}{t\log^{\star}t}+\frac{CC_1}{2}\dPC{}(\zeta,\sigma).\] 
      \qed
    \end{pr}

    We are now ready to control $N$ with the following lemma.

    \begin{lem}\label{boundN} There exists a constant $C_5$ such that
      \[N\leq C_5\log^{\star}\dsing{E}{x}Re^{C_4R},\]
    \end{lem}

    \begin{pr} By Lemma~\ref{bndRj+1-Rj}, if $j\in\llbracket0,N-2\rrbracket$, we have
      \[R_{j+1}-R_j\geq\frac{2\eps_0\dsing{E}{x_j}}{\underset{\left[\xi_j,\xi_{j+1}\right]}{\sup}\left(\eta\circ u_x\right)}.\]

      Denote by $\Sigma=\underset{\left[\xi_j,\xi_{j+1}\right]}{\sup}\left(\dsing{E}{u_x}\right)$. The function $t\mapsto t\log^{\star}t$ is increasing on $\intero{0}{1}$, thus
      \[R_{j+1}-R_j\geq\frac{2\eps_0C_1^{-1}\dsing{E}{x_j}}{\Sigma\log^{\star}\Sigma}.\]

      Since for $\zeta\in\left[\xi_j,\xi_{j+1}\right]$, $\dhermf{x}(u_x(\zeta),x_j)\leq\eps_0\dsing{E}{x_j}$, we have $\Sigma\leq(1+\eps_0)\dsing{E}{x_j}$. It follows that
      \[R_{j+1}-R_j\geq\frac{2\eps_0C_1^{-1}}{1+\eps_0}\frac{1}{\log^{\star}\Sigma}.\]

      Moreover, by Lemma~\ref{nottoocloseS}, $\log^{\star}\Sigma\leq e^{C_3+C_4R}\log^{\star}\dsing{E}{x}$. Hence,
      \[R_{j+1}-R_j\geq\frac{2\eps_0C_1^{-1}e^{-C_3-C_4R}}{(1+\eps_0)\log^{\star}\dsing{E}{x}}.\]

      This bound is uniform in $j$, and it is clear that
      \[N\leq\frac{R}{\inf_j(R_{j+1}-R_j)}\leq\frac{C_1e^{C_3}(1+\eps_0)}{2\eps_0}\log^{\star}\dsing{E}{x}Re^{C_4R}.\]
      \qed
    \end{pr}

    We are now able to construct the function $\psi$ that we need to apply Lemma \ref{Rdelimpsemilamb}.

    \begin{lem} \label{lemconstructionpsi} Suppose that $x,y\in\mani{M}\backslash\mani{E}$ are such that \eqref{hypmainthm} for a sufficiently large $R$. Then, there exists a smooth function $\psi\colon\adhDR{R}\to\leafu{y}$ without critical point, such that $\dhermset{\adhDR{R}}(\psi,u_x)\leq e^{-2R}$ and $\norm{d\psi}_{\infty}\leq 2A$. Here, we consider the norm of $d\psi$ with respect to the Poincar\'{e} metric on the source $\adhDR{R}$ and the restriction of $g_{\mani{M}}$ on the goal $\leafu{y}$.
    \end{lem}

    \begin{pr} We keep the notations that we introduced before Lemma \ref{bndRj+1-Rj}. The condition that $\dhermf{x}\left(u_x(\zeta),x_j\right)\leq\eps_0\dsing{E}{x_j}$ for all $\zeta\in\left[\xi_j,\xi_{j+1}\right]$ allows us to use Lemma~\ref{geomsetuporthproj}. This will ensure the existence of successive orthogonal projection $\Phi_{x_{j+1}y_{j+1}}\colon\ballleafu{x_{j+1}}{\eps_0\dsing{E}{x_{j+1}}}\to \ballleafu{y_{j+1}}{k\eps_0\dsing{E}{x_{j+1}}}$, with by induction $y_{j+1}=\Phi_{x_jy_j}(x_{j+1})$ and $y_0=\Phi_{xy}(x)$. We just have to check that at each step, $x_j$ and $y_j$ are close enough. As long as $\dherm{}(x_j,y_j)\leq\eps_1\dsing{E}{x_j}$, we have $\dherm{}(x_j,y_j)\leq e^{jK}\dherm{}(x,y_0)\leq e^{jK}\dherm{}(x,y)$. Hence
    \begin{equation}\label{majdxjyj}\frac{\dherm{}(x_j,y_j)}{\dsing{E}{x_j}}\leq e^{NK}\frac{\dherm{}(x,y)}{\dsing{E}{x_j}}\leq\exp\left(\log^{\star}\dsing{E}{x}\left(C_3+(KC_5R+1)e^{C_4R}-e^{\alpha^{-1}R}\right)\right),\end{equation}
    using our hypothesis~\eqref{hypmainthm}, Lemmas~\ref{nottoocloseS} and~\ref{boundN}. Let us suppose that $\alpha<C_4^{-1}$. If $R$ is sufficiently large, we have $\dherm{}(x_j,y_j)\leq\eps_1\dsing{E}{x_j}$ at each step. If follows that we can define $\psi\colon\left[0,\xi\right]\to\leafu{y}$ by gluing the successive $\Phi_{x_jy_j}\circ u_x$ on $\left[\xi_j,\xi_{j+1}\right]$. We can do the same work for each $\xi\in\set{D}$ such that $\dPC{}(0,\xi)=R$ and obtain a function $\psi\colon\adhDR{R}\to\leafu{y}$. We could imagine that this does not give a smooth function and not even a continuous function. Take $\xi^1,\xi^2\in\adhDR{R}$ with $\dPC{}(0,\xi^1)=\dPC{}(0,\xi^2)=R$. Define $\xi^1_j$, $\xi^2_j$ to have the above conditions and $\Cmod{\xi_j^1}=\Cmod{\xi_j^2}$ for both segments $\left[0,\xi^1\right]$ and $\left[0,\xi^2\right]$. If we suppose that $\Cmod{\xi^1-\xi^2}\leq e^{-R}A^{-1}\eps_0\exp\left(-e^{C_4R}\log^{\star}\dsing{E}{x}\right)$, denote by $x_j^k=u_x(\xi_j^k)$ and $y_j^k=\psi(\xi_j^k)$ for $j\in\llbracket0,N\rrbracket$ and $k\in\{1,2\}$. By definition, at first step $y_1^k=\Phi_{xy}(x_1^k)$, thus
    \[\dhermf{y}(y_1^1,y_1^2)\leq k\dhermf{x}(x_1^1,x_1^2)\leq kAe^{R}\Cmod{\xi_1^1-\xi_1^2}\leq kAe^R\Cmod{\xi^1-\xi^2}\leq k\eps_0\dsing{E}{x_1^1}.\]
    
    By Lemma~\ref{geomsetuporthproj}, it follows that $\Phi_{x_1^1y_1^1}=\Phi_{x_1^2y_1^2}$ on the intersection of their domains of definition. We continue by induction. By the same computation, we have at each step that
    \[\dhermf{y}(y_j^1,y_j^2)\leq k\eps_0\dsing{E}{x_j^1}.\]

    Hence, $\psi$ is locally equal to some $\Phi_{x_jy_j}\circ u_x$. Then, it is smooth. Moreover, by \eqref{majdxjyj}, if $R$ is sufficiently large,
    \begin{equation}\label{dgpsiux} \dhermset{\adhDR{R}}(\psi,u_x)\leq e^{NK}\dherm{}(x,y)\leq\exp\left(\log^{\star}\dsing{E}{x}\left(KC_5Re^{C_4R}-e^{\alpha^{-1}R}\right)\right)\leq e^{-2R}.\end{equation}

    Since $\psi$ is locally equal to $\Phi_{x_jy_j}\circ u_x$, it is clear that $\psi$ has no critical point. The condition on $\norm{d\psi}_{\infty}$ will be quite easy to obtain. Observe that locally, $\psi=\Phi_{x_jy_j}\circ u_x$. By~\eqref{majdxjyj} and Lemma~\ref{geomsetuporthproj}, $\norm{d\Phi_{x_jy_j}}_{\infty}\leq2$ if $R$ is sufficiently large. It follows that locally
    % \begin{equation}\label{majdpsi}\norm{d\psi}_{\infty}\leq\norm{d\Phi_{x_jy_j}}_{\infty}\norm{du_x}_{\infty}\leq 2A.\end{equation}
    \[\norm{d\psi}_{\infty}\leq\norm{d\Phi_{x_jy_j}}_{\infty}\norm{du_x}_{\infty}\leq 2A.\]
    \qed
    \end{pr}
    
    Denote by $\tilde{y}=\Phi_{xy}(x)$ and define $\Psi\colon\adhDR{R}\to\set{D}$ to be the lifting of $\psi$ \emph{via} $u_{\tilde{y}}$ satisfying $\Psi(0)=0$. By definition, $\psi=u_{\tilde{y}}\circ\Psi$. To satisfy the hypothesis of Lemma~\ref{Rdelimpsemilamb}, it remains to control $\norm{\mu_{\Psi}}_{\class{1}}$.

    \begin{lem}\label{lemmuPsi} With the hypothesis of Lemma \ref{lemconstructionpsi}, $\norm{\mu_{\Psi}}_{\class{1}}\leq e^{-2R}$.
    \end{lem}

    In order to prove this result, we will need first the following estimate on the second derivative of a uniformization.

        \begin{lem}\label{bounduzsec} There exist a radius $r_0\in\set{R}_+^*$ and a constant $C_6\in\set{R}_+^*$ such that for all $z\in\manis{M}{E}$, $u_z$ a uniformization of $\leafu{z}$ and $\Cmod{t}\leq\frac{r_0}{2}$, $\norm{u_z''(t)}\leq C_6$.
    \end{lem}

    \begin{pr} Note that $\dhermf{z}{}(z,u_z(t))\leq A\dPC{}(0,t)$. Hence, $\dPC{}(0,t)\leq R_0=\frac{d}{A}$ implies that $u_z(t)$ is contained in a holomorphic chart of $\mani{M}$. Consider $r_0=\frac{e^{R_0}-1}{e^{R_0}+1}$, $\Cmod{t}\leq\frac{r_0}{2}$ and apply Cauchy formula. We get
      \[\begin{aligned}\norm{u_x''(t)}_{g_{\mani{M}}}\leq C\norm{u_x''(t)}&\leq\frac{C}{\pi r_0}\int_0^{2\pi}\norm{u_x'\left(r_0e^{i\theta}\right)}d\theta\\ &\leq\frac{CC'(1-r_0^2)}{\pi r_0}\int_0^{2\pi}\frac{\norm{u_x'\left(r_0e^{i\theta}\right)}_{g_{\mani{M}}}}{1-r_0^2}d\theta\\ &\leq\frac{2CC'A(1-r_0^2)}{r_0},\end{aligned}\]
    where the constants $C$ and $C'$ come from the equivalence of the standard Hermitian metric and $g_{\mani{M}}$ in the holomorphic chart.
      \qed
    \end{pr}

\emph{Proof of Lemma \ref{lemmuPsi}.} Recall that locally, $\Psi=u_{\tilde{y}}^{-1}\circ\Phi_j\circ u_x$, where $\Phi_j=\Phi_{x_jy_j}$. An explicit computation and Lemma~\ref{geomsetuporthproj} give that
    \begin{equation}\label{boundmuPsiponct}\begin{aligned}\Cmod{\mu_{\Psi}}+\Cmod{\der{\mu_{\Psi}}{t}}+\Cmod{\bder{\mu_{\Psi}}{t}}&\leq\left(1+\frac{\norm{\left(\overline{\partial}\Phi_j\right)\circ u_x}}{\norm{\left(\partial\Phi_j\right)\circ u_x}}\right)\left(\frac{\norm{u_x'}}{\norm{\left(\partial\Phi_j\right)\circ u_x}}\norm{D^2\Phi_j\circ u_x}\right.\\ &\left.+2\frac{\norm{u_y''\circ\Psi}\norm{u_x'}}{\norm{u_y'\circ\Psi}^2}\norm{\left(\overline{\partial}\Phi_j\right)\circ u_x}\right)+\left(1+2\frac{\norm{u_x''}}{\norm{u_x'}}\right)\frac{\norm{\left(\overline{\partial}\Phi_j\right)\circ u_x}}{\norm{\left(\partial\Phi_j\right)\circ u_x}}
        \\&\leq3e^K\frac{\dhermset{\adhDR{R}}(u_x,\psi)}{\dsing{E}{x_j}}\left(1+\frac{\norm{u_x'}}{\dsing{E}{x_j}}+\frac{\norm{u_{\tilde{y}}''\circ\Psi}\norm{u_x'}}{\norm{u_{\tilde{y}}'\circ\Psi}^2}+\frac{\norm{u_x''}}{\norm{u_x'}}\right).\end{aligned}\end{equation}
For $t\in\set{D}$, consider the automorphism $f_t\colon z\mapsto\frac{z+t}{1+\cjg{t}z}$ of $\set{D}$, and define $v_z=u_z\circ f_t=u_{u_z(t)}$. We get
    \[u_z'(t)=\frac{v_z'(0)}{1-\Cmod{t}^2},\qquad u_z''(t)=\frac{v_z''(0)}{\left(1-\Cmod{t}^2\right)^2}+\frac{2\cjg{t}v_z'(0)}{\left(1-\Cmod{t}^2\right)^2}.\]

    Applying this to $z=x$ and $z=\tilde{y}$, using Brody hyperbolicity, Lemma~\ref{bounduzsec} and Proposition~\ref{setupbndlamb} in~\eqref{boundmuPsiponct}, we obtain
    \[\begin{aligned}\Cmod{\mu_{\Psi}(t)}+\Cmod{\der{\mu_{\Psi}}{t}(t)}+\Cmod{\bder{\mu_{\Psi}}{t}(t)}\leq \frac{3e^K}{1-\Cmod{t}^2}\frac{\dhermset{\adhDR{R}}(u_x,\psi)}{\dsing{E}{x_j}}&\left(1+\frac{A}{\dsing{E}{x_j}}+\frac{C_6A+2A^2}{C_1^2\left(\dsing{E}{\psi}\log^{\star}\dsing{E}{\psi}\right)^2}\right.\\ &+\left.\frac{C_6+2A}{C_1\dsing{E}{u_x}\log^{\star}\dsing{E}{u_x}}\right).\end{aligned}\]

    Using Lemma~\ref{nottoocloseS} on $\dsing{E}{x_j}$, $\dsing{E}{u_x}$ and $\dsing{E}{\psi}$, together with~\eqref{dgpsiux} and $1-\Cmod{t}^2\geq e^{-R}$, we obtain
    % \begin{equation}\label{boundnormC1muPsi}\norm{\mu_{\Psi}}_{\class{1}}\leq C\exp\left(\log^{\star}\dsing{E}{x}\left(C'Re^{C_4R}-e^{\alpha^{-1}R}\right)\right)\leq e^{-2R},\end{equation}
    \[\norm{\mu_{\Psi}}_{\class{1}}\leq C\exp\left(\log^{\star}\dsing{E}{x}\left(C'Re^{C_4R}-e^{\alpha^{-1}R}\right)\right)\leq e^{-2R},\]
    if $R$ is sufficently large.
    \qed

    \emph{End of proof of Theorem \ref{mainthm}.} Lemmas \ref{lemconstructionpsi} and \ref{lemmuPsi} ensure us that we are in the setup to use Lemma~\ref{Rdelimpsemilamb}. Since $x$ and $y$ play a symmetric role, if we denote by $\tilde{y}=\Phi_{xy}(x)$ and $\tilde{x}=\Phi_{yx}(y)$, we get
    \begin{equation}\label{lambdax'y}\eta(x)-\eta(\tilde{y})\leq C_2e^{-R},\qquad\eta(y)-\eta(\tilde{x})\leq C_2e^{-R}.\end{equation}

    It remains only to compare $\eta$ in two close points in the same leaf. That is, showing that $\eta(x)$ (resp. $\eta(y)$) and $\eta(\tilde{x})$ (resp. $\eta(\tilde{y})$) are close to each other. By symmetry of $x$ and $y$, we only show it for $x$ and $\tilde{x}$. Let $t\in\set{D}$ be a time such that $u_x(t)=\tilde{x}$ and $\dPC{}(0,t)\leq2\dPC{}(x,\tilde{x})$. By Lemma~\ref{geomsetuporthproj} and Proposition~\ref{setupbndlamb}, we have
    \begin{equation}\label{comptr0}\Cmod{t}\leq\dPC{}(0,t)\leq2\dPC{}(x,\tilde{x})\leq C\frac{\dhermf{x}(x,\tilde{x})}{\dsing{E}{x}\log^{\star}\dsing{E}{x}}\leq Ck\frac{\dherm{}(x,y)}{\dsing{E}{x}}\leq\frac{r_0}{2},\end{equation}
      if $R$ is sufficiently large. Hence, $x$ and $\tilde{x}$ belong to the same Hermitian chart and by Lemma~\ref{bounduzsec}, $\norm{u_x''}_{g_{\mani{M}}}\leq C$ on $\left[0,t\right]$. We have
      \[\begin{aligned}\Cmod{\eta(x)-\eta(\tilde{x})}&=\Cmod{\norm{u_x'(0)}_{g_{\mani{M}},x}-\frac{\norm{u_x'(t)}_{g_{\mani{M}},\tilde{x}}}{1-\Cmod{t}^2}}\\ &\leq\frac{4}{4-r_0^2}\Cmod{\norm{u_x'(0)}_{g_{\mani{M}},x}-\norm{u_x'(t)}_{g_{\mani{M}},\tilde{x}}}+\frac{4A\Cmod{t}^2}{4-r_0^2}\\
          &\leq CA\dherm{}(x,\tilde{x})+\frac{4A\Cmod{t}^2}{4-r_0^2}+\Cmod{\norm{u_x'(0)}_{g_{\mani{M}},\tilde{x}}-\norm{u_x'(t)}_{g_{\mani{M}},\tilde{x}}}\\
          &\leq C'\frac{\dherm{}(x,y)}{\dsing{E}{x}}+C''\norm{u_x'(0)-u_x'(t)}\leq C'\frac{\dherm{}(x,y)}{\dsing{E}{x}}+C''C_6\Cmod{t}\leq C\frac{\dherm{}(x,y)}{\dsing{E}{x}},\end{aligned}\]
      where we used successively the Cauchy formula (see just after \eqref{estimateCauchycomplambda}), the Lipschitz spatial variation of $g_{\mani{M}}$, the equivalence of Hermitian metrics, \eqref{comptr0}, Lemma \ref{geomsetuporthproj} and Lemma \ref{bounduzsec}. It follows that if $R$ is sufficiently large, $\Cmod{\eta(x)-\eta(\tilde{x})}\leq e^{-R}$. This, with the symmetric statement $\Cmod{\eta(y)-\eta(\tilde{y})}\leq e^{-R}$, together with~\eqref{lambdax'y}, give us the result of Theorem~\ref{mainthm}.
      \qed

      We finish our paper with an example that illustrates our enhancement to Theorem~\ref{thmDNSII}. It is a foliation on $\projective{2}$ that satisfies the hypothesis of Theroem~\ref{mainthm}, but at least one of its singularities is non-linearizable.  We leave the details of the easy computations to the reader.

      \begin{exmp} \label{exfin} Consider in $\set{C}^2$ the polynomial vector field
        \[X(z,w)=\left(2z+w^2-z^3\right)\der{}{z}+\left(w-z^2w\right)\der{}{w}.\]

      The field $X$ defines a singular holomorphic foliation by Riemann surfaces on $\set{C}^2$ that extends to a foliation $\fol{}=\left(\projective{2},\leafatlas,E\right)$ of degree 2 on $\projective{2}$. Computing the vector fields that define $\fol{}$ in the two other affine charts of $\projective{2}$, one can check that $\fol{}$ does not have any singularity on the line at the infinity. Moreover, solving $X=0$, one gets that
      \[E=\left\{(0,0);(\sqrt{2},0);(-\sqrt{2},0);(1,i);(1,-i);(-1,1);(-1,-1)\right\}.\]
      
      Computing the Jacobian matrix of $X$ in each of these singularities, one can check easily that every singularity of $\fol{}$ is non-degenerate. By the results of Glutsyuk~\cite{Glu} and Lins Neto~\cite{LN2}, it follows that $\fol{}$ is a Brody hyperbolic foliation on $\projective{2}$ and Theorem~\ref{mainthm} applies.
          
          We claim that $X$ is not linearizable near $0$. In fact, $X$ is essentially the most basic example of resonance (see~\cite{MR} or~\cite{Dulac} for more details about this phenomenon). In order to prove our claim, we blow up the origin and compare the forms of the liftings of $X$ and its linear part $X_{\ell}=2z\der{}{z}+w\der{}{z}$. Consider one of the usual projections $\pi(t,w)=(tw,w)$ from the total space of the blowing-up to $\set{C}^2$. On the one hand,
          \[\pi^*X_{\ell}=t\der{}{t}+w\der{}{w};\]
          and on the other hand,
          \[\pi^*X=(t+w)\der{}{t}+(w-t^2w^3)\der{}{w}.\]
          These two vector fields cannot define the same foliation. Indeed, their linear parts are not equivalent up to multiplication by $\lambda\in\set{C}^*$. Hence the singularity of $\fol{}$ at $0$ is not linearizable.

        \end{exmp}

\end{document}